\documentclass{amsart}
\usepackage{amssymb,hyperref}

\newcommand{\al}{\alpha}\newcommand{\be}{\beta}
\newcommand{\de}{\delta}
\newcommand{\ep}{\epsilon}\newcommand{\vep}{\varepsilon}

\newcommand{\ga}{\gamma}
\def\<{\langle}
\def\>{\rangle}
\newcommand{\R}{\mathbb{R}}\newcommand{\Z}{\mathbb{Z}}
\newcommand{\N}{\mathbb{N}}

\def\CO{\mathcal {O}}

\newcommand{\pt}{\partial_t}\newcommand{\pa}{\partial}

\newcommand{\les}{{\lesssim}}

\newcommand{\beeq}{\begin{equation}}\newcommand{\eneq}{\end{equation}}

\newtheorem{thm}{Theorem}[section]
\newtheorem{prop}[thm]{Proposition}
\newtheorem{coro}[thm]{Corollary}
\newtheorem{rem}{Remark}[section]

\newtheorem{lem}[thm]{Lemma}

\newenvironment{prf}{\noindent {\bf Proof.} }{\endprf\par}
\def \endprf{\hfill  {\vrule height6pt width6pt depth0pt}\medskip}

\numberwithin{equation}{section}

\begin{document}

\title[The Wave equation on Asymptotically Euclidean Manifolds]{Concerning the Wave equation on
Asymptotically Euclidean Manifolds}

\thanks{The first author was supported by the National Science Foundation.  The second author was supported in part by NSFC 10871175}

\author{Christopher D. Sogge}
\address{Department of Mathematics, Johns Hopkins University, Baltimore,
Maryland 21218}
\email{sogge@jhu.edu}


\author{Chengbo Wang}
\address{Department of Mathematics, Johns Hopkins University, Baltimore,
Maryland 21218} \email{wangcbo@jhu.edu}

\dedicatory{} \commby{}

\begin{abstract}
We obtain KSS, Strichartz and certain weighted Strichartz estimates for the wave equation on $( \R^d , \mathfrak{g})$, $d \geq 3$, when metric $\mathfrak{g}$ is non-trapping and approaches the Euclidean metric like $\< x \>^{- \rho}$ with $\rho>0$. Using the KSS estimate, we prove almost global existence for quadratically semilinear wave equations with small initial data for $\rho> 1$ and $d=3$.
Also, we establish the Strauss conjecture when the metric is radial with $\rho>0$ for $d= 3$.
\end{abstract}

\maketitle\tableofcontents

\section{Introduction and Main Results}\label{7-sec-intro}
This paper is devoted to the study of the semilinear wave equation on asymptotically Euclidean non-trapping Riemannian manifolds. We shall obtain almost global existence for quadratic semilinear wave equations with small data, and show that the
Strauss conjecture holds in this setting, in dimension $d=3$.

In Minkowski space, the quadratically semilinear wave equation has
been thoroughly studied. Global existence is known in dimension
$d\ge 4$ for small initial data (see Klainerman and Ponce
\cite{KlPo83_01} and references therein). Almost global existence in
dimension $d=3$ for small data was shown by John and Klainerman in
\cite{JoKl84_01}. Almost global means that the life time of a
solution is at least $\exp(c/ \delta)$ with some $c>0$, where
$\delta$ is the size of the initial data in some appropriate Sobolev
space. Note that, in dimension $d=3$, Sideris~\cite{Si83_01} has
proved that global existence does not hold in general (see also
John~\cite{Jo81_01}).

In \cite{KeSmSo02}, Keel, Smith and Sogge gave a new proof of the
almost global existence result in dimension $3$ using estimates
(known as KSS estimates) of the form
\begin{equation} \label{70-est-KSSE}
( \ln (2+ T ) )^{-1/2} \| \< x \>^{-1/2} u' \|_{L^2 ( [0,T] \times \R^3 )}
\lesssim \| u' (0, \cdot ) \|_{L^2 ( \R^3 )} + \int_0^T \| F(s, \cdot ) \|_{L^2 ( \R^3 )} d s ,
\end{equation}
and a certain Sobolev type estimate due to Klainerman (see
\cite{Kl85_01}). Here $u$ solves the wave equation $\Box u = F$ in
$[0 , + \infty ) \times \R^d$ and $u' = ( \partial_{t} u ,
\partial_{x} u )$.
Earlier versions of \eqref{70-est-KSSE} appeared before.  The first ones appear
to be due to Morawetz \cite{Mor68} and Strauss~\cite{Strauss},
who proved  somewhat weaker versions of
\eqref{70-est-KSSE}.  See also \cite{KPV}.
In \cite{KeSmSo02} existence results for the non-trapping obstacle case were
also obtained.
In \cite{KeSmSo04_01}, similar results were obtained for the
corresponding quasilinear equation (see also Metcalfe-Sogge
\cite{MeSo06_01}).

Recently, Bony and H\"{a}fner \cite{BoHa} obtained a weaker version
of KSS estimate in the current setting and proved the long time existence
for quadratic semilinear wave equations with small data. In the
present  paper, by using results of Metcalfe-Sogge \cite{MeSo06_01},
we are able to prove the full KSS estimate, and hence the almost
global existence for the quadratic semilinear wave equation.

Recently, in Minkowski space, Fang and Wang \cite{FaWa} and
Hidano-Metcalfe-Smith-Sogge-Zhou \cite{HMSSZ} proved the Strauss
conjecture with low regularity for $d=2,3,4$, by using a weighted
Strichartz estimate of the form
\beeq\label{70-est-weighted-Strichartz}
\||x|^{\frac{n}{2}-\frac{n+1}{r}-\ga} u \|_{L^r_{t, |x|} L^2_\omega
} \les \|u(0, \cdot)\|_{\dot{H}^{\ga}_x}+\|\pt
u(0,\cdot)\|_{\dot{H}^{\ga-1}_x}+ \| F \|_{L^1_t\dot{H}^{\ga-1}_x}
\eneq
for $\ga\in (1/2-1/r,d/2-1/r) $, $r\in [2,\infty]$, where we use the notation
$$\|f\|_{L^q_t L^r_{|x|} L^p_\omega}= \|(\int_0^\infty \|f(t,|x|\omega)\|^r_{L^p_\omega} |x|^{d-1} d |x|)^{1/r}\|_{L^q_t},$$
with $L^p_\omega$ denoting the $L^p$-norm on $S^{n-1}$ with respect to the standard
measure.
In this paper, we obtain a somewhat weaker version of this estimate
in this general setting, which suffices for us to establish the
Strauss conjecture for $d=3$ (when the metric is radial).

Using ideas from Burq \cite{Burq}, Metcalfe \cite{Metcalfe},
Smith-Sogge  \cite{SmSo00} and Hidano-Metcalfe-Smith-Sogge-Zhou
\cite{HMSSZ}, we can also use the local energy decay estimates to
prove global Strichartz estimates in this setting.  We should point
out, though, that the idea that, in many situations, local energy
estimates can be used to prove global Strichartz estimates occurs in
many other works.  The first seems to be that of Journ\'e, Soffer
and Sogge \cite{JSS} who proved global Strichartz estimates for
Schr\"odinger operators with potential using local energy estimates
(local smoothing) for $e^{it\Delta}$.  Staffilani and Tataru
\cite{StTa} extended this philosphy by considering more general
perturbations of $\Delta$, and more recently Metcalfe and Tataru
\cite{MeTa07p} used the philosophy that local energy estimates imply
Strichartz estimates to handle (small) metric perturbations of
$\square=\partial_t^2-\Delta$.  Thus, in many ways, some of the
techniques employed in this paper are not novel, since they have
been used in many earlier works.  A slight novelty, though, might be
that we obtain our global estimates by combining local energy
estimates (in this case due to Bony and H\"afner \cite{BoHa}) with
global Strichartz estimates not involving $\Delta$, but rather small
perturbations of the Laplacian.  The ones that allow us to prove the
aforementioned (sharp) KSS estimates are due to Metcalfe and the
first author \cite{MeSo06_01}, while the ones that allow us to prove
the standard mixed-norm Strichartz estimates are due to Metcalfe and
Tataru~\cite{MeTa07p}.

Let us now state our precise results. We consider asymptotically
Euclidean manifolds $( \R^d , \mathfrak{g})$ with $d \geq 3$ and
\begin{equation*}
\mathfrak{g} = \sum_{i,j=1}^{d} g_{ij} (x) \, d x^i \, d x^j .
\end{equation*}
We suppose $g_{ij} (x) \in C^{\infty} ( \R^{d} )$ and, for some $\rho >0$,
\begin{equation}\tag{H1} \label{c1}
\forall \alpha \in \N^d \qquad \partial^{\alpha}_x ( g_{ij} -
\delta_{ij} ) = \CO ( \< x \>^{- \vert \alpha \vert - \rho} ) ,
\end{equation}
with $\delta_{ij}=\delta^{ij}$ being the Kronecker delta function.
We also assume that
\begin{equation}\tag{H2} \label{c15}
\mathfrak{g} \text{ is non-trapping.}
\end{equation}
Let $g (x) = ( \det ( \mathfrak{g} ) )^{1/4}$. The Laplace--Beltrami
operator associated with $\mathfrak{g}$ is given by
\begin{equation*}
\Delta_{\mathfrak{g}} = \sum_{ij} \frac{1}{g^2} \partial_i g^{ij} g^2 \partial_j ,
\end{equation*}
where $g^{ij} (x)$ denotes the inverse metric. Note
$-\Delta_{\mathfrak{g}}$ is self-adjoint non-negative on $L^2(\R^d,
g^2 dx)$, while $P=-g\Delta_{\mathfrak{g}}g^{-1}$ is self-adjoint
non-negative on $L^2(\R^d, dx)$. Let $\Omega = \Omega_{k , \ell} : =
x_{k} \partial_{\ell} - x_{\ell} \partial_{k}$ be the rotational
vector fields. We consider the following semilinear wave equation
\begin{equation} \label{70-eqn-DSLW}
\begin{cases}
\Box_{\mathfrak{g}} u=Q ( u'),
\quad (t,x)\in \R_+\times \R^d
\\
u(0,x)=u_0(x), \quad \partial_t u(0,x)=u_1(x), \quad x\in \R^d.
\end{cases}
\end{equation}
Here $\Box_{\mathfrak{g}} = \partial_{t}^{2}-\Delta_{\mathfrak{g}}$
and $Q (u')$ is a quadratic form in $u' = ( \partial_t u ,
\partial_{x} u )$. One of our main results is the following theorem.

\begin{thm}\sl  \label{70-thm-DSLW-AGlob}
Assume hypotheses \eqref{c1} and \eqref{c15} with $\rho>1$. Suppose
$u_0 , u_1 \in C_0^{\infty} ( \R^3 )$, and
\begin{equation}  \label{70-eqn-DSLW-data}
\sum_{\vert \alpha \vert + |\beta| \leq 5} \big\Vert \partial^\be_x
\Omega^{\alpha} u_0 \big\Vert_{L^2_x} + \sum_{\vert \alpha \vert +
|\be| \leq 4} \big\Vert \partial^\be_x \Omega^{\alpha} u_1 \big\Vert_{L^2_x} \leq \delta .
\end{equation}
For $\delta$ small enough, the problem \eqref{70-eqn-DSLW} has a
unique almost global solution $u \in C^{\infty} ( [ 0 , T_\de]
\times \R^d )$ with $T_\delta=exp(c/\delta)$ for some $c>0$.
\end{thm}

The main ingredient of the proof are estimates of type
\eqref{70-est-KSSE}. Let us therefore consider the corresponding
linear equation. Let $u$ be solution of
\begin{equation} \label{70-eqn-LW}
\begin{cases}
(\partial_t^2-\Delta_{\mathfrak{g}})u(t,x)=F(t,x),
\quad (t,x)\in \R_+\times \R^d
\\
u(0,x)=u_0(x), \quad \partial_t u(0,x)=u_1(x), \quad x\in \R^d.
\end{cases}
\end{equation}
With the notation
\begin{equation*}
A_{\mu} (T) = \left\{ \begin{aligned}
&\ln (2+T)^{-1/2} &&\mu = 1/2 ,   \\
&1 &&\mu > 1/2 ,
\end{aligned} \right.
\end{equation*}
we have the following KSS estimate.

\begin{thm}\sl \label{70-thm-KSS}
Assume that \eqref{c1} and \eqref{c15} hold with $\rho > 0$ and let
$ \mu\ge 1/2$, $d\ge 3$. For all $\varepsilon >0$, the solution of
\eqref{70-eqn-LW} satisfies
\begin{equation} \label{70-est-KSS}
A_{\mu}(T)\left(\big\Vert \< x \>^{-\mu} u' \big\Vert_{L^2_T
L^2_x}+\big\Vert \< x \>^{-\mu-1} u \big\Vert_{L^2_T L^2_x} \right)
\lesssim    \Vert u' ( 0 , \cdot ) \Vert_{L^2_x} +  \Vert \< x
\>^{\frac{1}{2}+\vep} F\Vert_{L^2_T L^2_x},
\end{equation} where we use $L^2_T$ to denote $L^2_{t\in [0, T]}$.
\end{thm}

We remark that our estimate \eqref{70-est-KSS} agrees with ones in
Bony and H\"afner \cite{BoHa} when $\mu>1/2$, while for $\mu=1/2$
they are slightly stronger since we obtain the sharp bounds with
$A_{1/2}(T)=(\ln (2+T))^{-1/2}$, as opposed to the bounds of
$T^{-\varepsilon}$, $\varepsilon>0$, in \cite{BoHa}.  This
improvement allows us to obtain the almost global existence results
alluded to before.  On the other hand, our proof is very similar to
that of Bony and H\"afner \cite{BoHa} and papers that preceded it,
starting with \cite{KeSmSo02}.  A slight point of departure is that
we combine local energy estimates (due to Bony and H\"afner
\cite{BoHa}) not with global KSS estimates for $\Delta$ but rather
for small metric perturbations of $\Delta$ (which are due to
Metcalfe and the first author \cite{MeSo06_01}).

To prove the nonlinear theorem, we need to get higher order
estimates.  For this purpose, let us put $Z=\{\partial_t ,
\partial_x, \Omega\}$, $Y= \{ \partial_x, \Omega\}$, $X = \{
\partial_x \}$. Then, we have

\begin{thm}\sl \label{70-thm-KSS-HighOrder}
Assume that \eqref{c1} and \eqref{c15} hold with $\rho > 1$.  Let $N
\ge 0$ and $ \mu\ge 1/2$. The solution of \eqref{70-eqn-LW}
satisfies
\begin{align}\label{70-est-KSS-HighOrder}
\sup_{0 \leq t \leq T} \sum_{1 \leq k + j \leq N+1} \big\Vert \partial_{t}^{k}
P^{j/2}g u (t & , \cdot ) \big\Vert_{L^2_x} + \sum_{\vert \alpha \vert \leq N}
A_{\mu}(T) \big\Vert \<x\>^{-\mu} (|(Z^{\alpha} u)'|+  \frac{|Z^{\alpha} u|}{\<x\>})
 \big\Vert_{L^2_T L^2_x} \nonumber \\
& \lesssim \sum_{\vert \alpha \vert \leq N} \big\Vert (Y^{\alpha}u)' (0, \cdot ) \big\Vert_{L^2_x}
 + \sum_{\vert \alpha \vert \leq N}\int_0^T \big\Vert Z^{\alpha} F(s, \cdot ) \big\Vert_{L^2_x} d s .
\end{align}
\end{thm}

Note that the estimate \eqref{70-est-KSS} can be viewed as the local
energy decay estimate for  $\Box_{\mathfrak{g}}$. The local
Strichartz estimates for variable coefficient wave equations have
been studied extensively, see e.g. Kapitanski \cite{Ka91},
Mockenhaupt-Seeger-Sogge \cite{MoSeSo93}, Smith \cite{Sm98},
Bahouri-Chemin \cite{BaCh99-1} \cite{BaCh99-2}, Tataru \cite{Ta00}
\cite{Ta01} \cite{Ta02}. And recently, Metcalfe and Tataru have
obtained global Strichartz estimates involving small perturbations of
the Minkowski metric in \cite{MeTa07p}.
As we mentioned before, by combining these with local energy estimates we
shall prove
global Strichartz estimates for
$\Box_{\mathfrak{g}}$.

For the Minkowski case, it is known (\cite{KeTa98}, \cite{LdSo95})
that we have Strichartz estimates if $(s,q,r)$ is admissible, i.e.,
$$\frac{1}{q}\le \min(\frac{1}{2}, \frac{d-1}{2}\left(\frac{1}{2}-
\frac{1}{r}\right)), \quad (q,r)\neq (2,\infty),
(\infty,\infty),\quad
s=d\left(\frac{1}{2}-\frac{1}{r}\right)-\frac{1}{q}.$$ Our global
Strichartz estimates for $\Box_{\mathfrak{g}}$ are the
following.

\begin{thm}[Global Strichartz estimate]\label{70-thm-Strichartz-Glob}
Assume that \eqref{c1} and \eqref{c15} hold with $\rho > 0$, $d\ge
3$ and $s\in [0,1]$ ($s\in (0,1)$ if $d=3$). The solution of
\eqref{70-eqn-LW} satisfies
\begin{equation}\label{70-est-Strichartz}
\big\Vert u  \big\Vert_{L^q_t L^r_x}\lesssim \|u_0\|_{\dot H^s}+
\|u_1\|_{\dot H^{s-1}}+\|F\|_{L^1_t \dot H^{s-1}},
\end{equation} for any admissible $(s,q,r)$ with $q>2$, $r<\infty$.
\end{thm}

The proofs of Theorem \ref{70-thm-KSS} and Theorem
\ref{70-thm-Strichartz-Glob} follow a similar strategy.  We first
use results from \cite{MeSo06_01} and \cite{MeTa07p} to show that we
can construct a metric $\tilde{\mathfrak{g}}$ which agrees with
$\mathfrak{g}$ near infinity and has the property that the bounds in
these two theorems are valid if $\square_\mathfrak{g}$ is replaced
by $\square_{\tilde{\mathfrak{g}}}$.  Then, by adapting arguments
from \cite{Burq} and \cite{SmSo00}, we can use these estimates along
with the local energy decay estimates for $\square_\mathfrak{g}$
(see Lemma \ref{70-thm-LocEnerDecay} below) to show that
$\square_\mathfrak{g}$ satisfies the same global estimates as its compact
perturbation $\square_{\tilde{\mathfrak{g}}}$.

Now let us describe the weighted Strichartz estimate and its
application to Strauss conjecture in this general setting. Let
$p>1$,
$$s_c=\frac{d}{2}-\frac{2}{p-1},\quad
s_{sb}=\frac{1}{2}-\frac{1}{p}.$$ The equation that we shall
consider is
\begin{equation}\label{70-eqn-SLW}
\begin{cases}
(\partial_t^2-\Delta_{\mathfrak{g}})u(t,x)=F_p(u(t,x)),
\quad (t,x)\in \R_+\times \R^d
\\
u(0,x)=u_0(x), \quad \partial_t u(0,x)=u_1(x), \quad x\in \R^d,
\end{cases}
\end{equation}
  We shall assume that the nonlinear term behaves like
$|u|^p$, 
and so we assume that
\begin{equation}\label{70-eqn-SLW-Fp}
\sum_{0\le j\le 1} |u|^j\, |\, \partial^j_u F_p(u)\, | \,
\lesssim \, |u|^p.
\end{equation}
See \cite{Sog95}, \S 4.4 for a discussion about how $s\ge s_{sb}$
is needed for local existence, while $s_c$ is critical for global existence.

We can now state our existence theorem for \eqref{70-eqn-SLW}. Due
to some technical difficulties, we are only able to deal with only the case where
$g^{ij}(x)=h(|x|)\delta^{ij}$ for some function $h$.

\begin{thm}\label{70-thm-Strauss} Assume that $g^{ij}(x)=h(|x|)\delta^{ij}$
for some function $h$, \eqref{c1} and \eqref{c15} hold with $\rho > 0$, $d=3$ and $p > p_c =1+\sqrt{2}$.
Then for any $\ep>0$ such that
\begin{equation}\label{70-eqn-Strauss-Regu}
s=s_c-\ep \in (s_{sb}, 1/2),
\end{equation}
there is an $\delta>0$ depending on $p$ so that
\eqref{70-eqn-SLW} has a global solution satisfying
$(Y^\alpha u(t,\cdot), \partial_t Y^\alpha u(t,\cdot))\in
\dot{H}^s \times \dot{H}^{s-1}$, $|\alpha|\le 1$, $t\in \R_+$,
whenever the initial data satisfies
\begin{equation}\label{70-eqn-SLW-data}
\sum_{|\alpha|\le 1}\left(\, \|Y^\alpha u_0\|_{\dot H^s} +\|Y^\alpha
u_1\|_{\dot H^{s-1}}\, \right)<\delta
\end{equation}
with $0<\delta<\delta_0$.  \end{thm}

Existence results of this type when $\Delta_g=\Delta$ are a celebrated result
of John \cite{Jo81_01}.  Subsequently, Strauss conjectured that for dimensions $d\ge2$
the critical exponent for small data global existence for equations of the form
\eqref{70-eqn-SLW} (when $\Delta_g=\Delta$) should be the positive root of the
equation $(d-1)p^2-(d+1)p-2=0$.  This conjecture was settled for the Minkowski
space case in  \cite{GLS97}, \cite{Glassey},  \cite{LS2}, \cite{Sideris} and \cite{Zhou}.
See \cite{Sog95} for further discussion.

As in the case of \eqref{70-eqn-DSLW}, the main ingredient of the
proof are estimates of type \eqref{70-est-weighted-Strichartz}.
If we consider the corresponding linear equation \eqref{70-eqn-LW}, then we
have the following estimate, where the metric is not restricted to
the special case where $g^{ij}(x)=h(|x|)\delta^{ij}$.

\begin{thm}\sl \label{70-thm-WSE}
Assume that \eqref{c1} and \eqref{c15} hold with $\rho > 0$, $d\ge
3$, $2< q\le \infty$ and $s\in (s_{sb}(q),1]$. For all $\varepsilon,
\eta >0$ small enough, the solution of \eqref{70-eqn-LW} satisfies
\begin{equation} \label{70-est-WSE}
\||x|^{\frac{d}{2}-\frac{d+1}{q}-s-\varepsilon} u
\|_{L^q_{t, |x|\ge R } L^{2+\eta}_\omega }+ \|\<x\>^{-\frac{1}{2}-s-\vep} u
\|_{L^2_{t, x}}\les
\|u_0\|_{\dot{H}^{s}_x}+\|u_1\|_{\dot{H}^{s-1}_x}+ \| F \|_{L^1_t\dot{H}^{s-1}_x}
\end{equation}
\end{thm}

We should comment on the hypotheses in the existence theorems.
First, because of the various commutator terms that arise in the proofs  we are, at present,
only able to handle semilinear terms in the existence theorems involving quadratic nonlinearities,
as opposed the the quasilinear case (for  $\square_g=\partial_t^2-\Delta$) treated in
\cite{KeSmSo04_01}.  For similar reasons, in our results involving the Strauss conjecture,
due to difficulties in dealing with commutators involving the $\Omega$ vector fields and $\square_\mathfrak{g}$
 we have to assume that the metric $\mathfrak{g}$ is spherically
symmetric.  For similar reasons, although the linear estimates just require the hypothesis
that $\rho>0$ currently our techniques require the assumption that $\rho>1$ in the
hypotheses of some of the nonlinear theorems.  We do not know, however, what the natural assumption
regarding $\rho$ should be for the latter, though.

\section{KSS Estimates}\label{70-sec-Kss}
In this section, we give the proof of the KSS estimates. First, we
will need the following lemmas, where we denote
$\widetilde{\partial}_x : = \partial_x g^{-1}$.
\begin{lem}[Theorem 1.3 and Proposition 4.6 in \cite{BoHa}]\label{70-thm-Morawetz}
Assume that \eqref{c1} and \eqref{c15} hold with $\rho > 0$, then
for all $\ep >0$, the solution of the equation $(\pt^2+P)u=F$
satisfies \beeq\label{70-est-Morawetz} \big\Vert \< x
\>^{-\frac{1}{2}-\ep} (\pt, P^{1/2}) u \big\Vert_{L^2 ( \R \times
\R^{d} )} \lesssim    \Vert (\pt, P^{1/2}) u ( 0 , \cdot )
\Vert_{L^2 ( \R^{d} )} +  \Vert \< x \>^{\frac{1}{2}+\ep}
F\Vert_{L^2 ( \R \times \R^{d} )}.\eneq
\end{lem}

\begin{rem}\label{70-rem-Morawetz}
  In fact, from the proof of Proposition 4.4 and 4.6 in \cite{BoHa}, we also have
   \beeq\label{70-est-Morawetz-LowOrder} \big\Vert \< x \>^{-\frac{3}{2}-\ep} u
  \big\Vert_{L^2 ( \R \times \R^{d} )} \lesssim    \Vert (\pt, P^{1/2})u ( 0 , \cdot ) \Vert_{L^2 ( \R^{d} )}
   +  \Vert \< x \>^{\frac{1}{2}+\ep} F\Vert_{L^2 ( \R \times \R^{d} )}.\eneq
\end{rem}

\begin{lem}[Theorem 5.1 in \cite{MeSo06_01}]\label{70-thm-KSS-SmallPert}
Let $\Box_h=\pt^2-\Delta+h^{\al\be}(t,x)\pa_\al \pa_\be$, $h^{\al\be}=h^{\be\al}$ and $\sum |h^{\al\be}|\le \de$.
Then if $\de>0$ is small enough, $d\ge 3$, the solution to the equation $\Box_h u=F$ satisfies
\begin{align}
(\ln(2+T))^{-1}\big\Vert \< x \>^{-1/2} (|u'|+\frac{|u|}{\<x\>}) \big\Vert^2_{L^2 ( [0, T] \times \R^{d} )}+
\big\Vert \< x \>^{-1/2-\ep} (|u'|+\frac{|u|}{\<x\>}) \big\Vert_{L^2 ( [0, T] \times \R^{d} )}^2
\nonumber \\
 \lesssim    \Vert u' ( 0 , \cdot ) \Vert^2_{L^2 ( \R^{d} )} +  \int^T_0 \int \left(|u'|+\frac{u}{|x|}\right)
 \left(|F| +(|h'|+\frac{h}{|x|})|u'|\right)d x dt\label{70-est-KSS-SmallPert}
\end{align} for any $\ep>0$.
\end{lem}

To obtain higher order estimates, we will need the following.
\begin{lem}[Lemma B.13, 4.1 and 4.2 in \cite{BoHa}]\label{70-lem-equi-PLapBH}
  For all $-3/2\le\widetilde{\mu} < \mu \leq 3/2$, we have
\beeq\label{a}\big\Vert \<x\>^{-\mu} \widetilde{\partial}_{\ell} u
\big\Vert_{L^2 ( \R^{d} )} \lesssim \big\Vert \<x\>^{-
\widetilde{\mu}} {P^{1/2}} u \big\Vert_{L^2 ( \R^{d} )} , \eneq
\beeq\label{b}\big\Vert \<x\>^{- \mu} {P^{1/2}} u \big\Vert_{L^2 (
\R^{d} )} \lesssim \sum_{\ell =1}^{d} \big\Vert \<x\>^{-
\widetilde{\mu}} \widetilde{\partial}_{\ell} u \big\Vert_{L^2 (
\R^{d} )} .\eneq Also, for $u \in H^{1} (\R^{d})$,
\begin{equation}\label{c}
\Vert {P^{1/2}} u \Vert_{L^2 ( \R^{d} )} \lesssim \Vert \nabla
g^{-1} u \Vert_{L^2 ( \R^{d} )} \lesssim \Vert {P^{1/2}}  u
\Vert_{L^2 ( \R^{d} )} .
\end{equation}
\end{lem}

We also need a lemma which says that the homogeneous spaces defined
by $P$ and $-\Delta$ are essentially  the same. In what follows,
``remainder terms", $r_{j}$, $j\in \N$, will denote a smooth
function such that \begin{equation} \label{c17}
\partial^{\alpha}_{x} r_{j} (x) = \CO \big( \< x \>^{-\rho - j - \vert \alpha \vert} \big) .
\end{equation}
\begin{lem}\label{70-lem-equi-PLap}
If $s\in[-1,1]$, then \beeq\label{d}\|u\|_{\dot{H}^s}\simeq
\|P^{s/2}u\|_{L^2_x}.\eneq If $s\in[0,1]$,
\beeq\label{f}\|\tilde{\pa}_j u\|_{\dot{H}^{-s}} \les \|P^{1/2}
u\|_{\dot{H}^{-s}},\eneq \beeq\label{e}\|P^{1/2} u\|_{\dot{H}^s}\les
\sum_j\|\tilde{\pa}_j u\|_{\dot{H}^s}.\eneq Moreover, we have for
$s\in (0,2]$ and $1<q <d/s$, \beeq\label{g}\|P^{s/2}u\|_{L^q_x}\les
\|u\|_{\dot{H}^{s,q}}. \eneq
\end{lem}
\begin{prf}
For the first estimate \eqref{d}, by interpolation and duality, we
need only to prove the estimate for the special case where $s=1$, i.e.
$$\|\nabla u\|_{L^2}\simeq \|\nabla g^{-1} u\|_{L^2}.$$
In fact, $$\|\nabla g^{-1}u\|_{L^2}\les\|(\nabla g^{-1})u\|_{L^2}+\|g^{-1}\nabla u\|_{L^2}\les( \|\nabla g^{-1}\|_{L^d}+\|g^{-1}\|_{L^\infty})\|\nabla u\|_{L^2}.$$
By the hypotheses \eqref{c1} and the ellipticity of $P$, we know that $$Image(g)\subset (\delta,\delta^{-1}),\quad |\pa g|, |\pa g^{-1}|=\mathcal{O}( \<x\>^{-1-\rho})\in L^d,\quad |\pa^2 g^{-1}|=\mathcal{O}( \<x\>^{-2-\rho})\in L^{d/2},$$ for some $\delta>0$.
So we know that $$\|\nabla h u\|_{L^2}\les \|\nabla u\|_{L^2}$$ whenever $h=g$ or $h=g^{-1}$.

To prove \eqref{f}, we note first that by the first inequality
\eqref{d},
$$\|\tilde{\pa}_j f\|_{\dot H^{-1}}=\|\pa_j g^{-1} f\|_{\dot
H^{-1}}\les \|g^{-1} f\|_{L^2}\les \| f\|_{L^2}\les \|P^{1/2}
f\|_{\dot H^{-1}}.$$ Thus the inequality \eqref{f} follows from
interpolation with \eqref{c}.

Note that by \eqref{f}, $$|\<P^{1/2}f,P^{1/2}h\>|=|\<P f,
h\>|=|\<g^2 g^{ij}\tilde{\pa}_j f, \tilde{\pa}_i h\>|\les
\sum_j\|\tilde{\pa}_j f\|_{\dot{H}^s} \|P^{1/2} h\|_{\dot{H}^{-s}}
,$$ this gives \eqref{e}.

For the last inequality, let $a^{ij}=g^2 g^{ij}$ and $P_1=\pa_i a^{ij} \pa_j$, then $P g=g^{-1}P_1$, $a^{ij}\in L^\infty$ and $\pa_i a^{ij}\in L^d$. Denoting $D=\sqrt{-\Delta}$, then for $1<q<d$
$$\|P_1 D^{-2}u\|_{L^q}\les \|(\pa_i a^{ij}) \pa_j D^{-2} u\|_{L^q}+\|a^{ij}D^{-2}\pa_i \pa_j u\|_{L^q}\les \|u\|_{L^q}.$$ Thus if $1<q<d/2$,
$$\|P u\|_{L^q}=\|g^{-1}P_1 g^{-1} u\|_{L^q}\les\|P_1 g^{-1}u\|_{L^q}\les  \|D^2 g^{-1} u\|_{L^q}\les \|u\|_{\dot{H}^{2,q}}.$$
Consequently, the last inequality in the lemma follows from interpolating with the trivial estimate where $s=0$.
\end{prf}

Using this we can obtain analogues of Lemma \ref{70-thm-Morawetz} and Remark \ref{70-rem-Morawetz} involving $-\Delta_\mathfrak{g}$.
\begin{coro}\label{70-thm-Morawetz2}
Assume that \eqref{c1} and \eqref{c15} hold with $\rho > 0$, then
for all $\ep >0$, the solution of the equation
$(\pt^2-\Delta_{\mathfrak{g}})u=F$ satisfies
\beeq\label{70-est-Morawetz2} \big\Vert \< x \>^{-\frac{1}{2}-\ep}
\pa_{t,x} u \big\Vert_{L^2_{t,x}} +\big\Vert \< x
\>^{-\frac{3}{2}-\ep} u \big\Vert_{L^2_{t,x}} \lesssim    \Vert
\pa_{t,x} u ( 0 , \cdot ) \Vert_{L^2_x} +  \Vert \< x
\>^{\frac{1}{2}+\ep} F\Vert_{L^2_{t,x}}.\eneq
\end{coro}
\begin{prf}
  Let $v=g u$ and $G=g F$, then we have $$(\pt^2-g\Delta_{\mathfrak{g}}g^{-1})v=G.$$
Then by Lemma \ref{70-thm-Morawetz}, Remark \ref{70-rem-Morawetz},
\eqref{a} and \eqref{c} in Lemma \ref{70-lem-equi-PLapBH}, we have
\begin{eqnarray*}
 \|\< x \>^{-\frac{1}{2}-\ep} \pa_{t,x} u \|_{L^2_{t,x}} +\| \< x
\>^{-\frac{3}{2}-\ep} u \|_{L^2_{t,x}} &\les&  \|\< x
\>^{-\frac{1}{2}-\ep} (\pt, \tilde\pa_{x}) v \|_{L^2_{t,x}} +\| \< x
\>^{-\frac{3}{2}-\ep} v \|_{L^2_{t,x}} \\
&\les& \|\< x \>^{-\frac{1}{2}-\ep/2} (\pt, P^{1/2}) v \|_{L^2_{t,x}} +\|
\< x \>^{-\frac{3}{2}-\ep/2} v \|_{L^2_{t,x}}\\
&\les& \| (\pt, P^{1/2}) v ( 0 , \cdot ) \Vert_{L^2_x} + \Vert \< x
\>^{\frac{1}{2}+\ep/2} G\Vert_{L^2_{t,x}}\\
&\les& \| (\pt, \tilde{\pa}_x) v ( 0 , \cdot ) \Vert_{L^2_x} + \Vert
\< x \>^{\frac{1}{2}+\ep} F\Vert_{L^2_{t,x}}\\
&\les& \| (\pt, \pa_x) u ( 0 , \cdot ) \Vert_{L^2_x} + \Vert \< x
\>^{\frac{1}{2}+\ep} F\Vert_{L^2_{t,x}}.\end{eqnarray*} This
completes the proof.
\end{prf}

We now establish the local energy decay estimates for the $\Box_{{\mathfrak{g}}}$.
\begin{lem}[Local Energy Decay]\label{70-thm-LocEnerDecay}
For the linear equation \eqref{70-eqn-LW}, if
$F(t,x)= 0$ for $|x|>R$ with $R$ fixed, then for fixed
$\beta\in C^\infty_0( \R^d)$, we have
\beeq\label{70-est-LocEnerDecay1} \big\Vert \beta u
\big\Vert_{L^2_{t}H^1}\lesssim    \Vert u_0\Vert_{\dot{H}^{1}}+\Vert
u_1\Vert_{L^2_x}+  \Vert F\Vert_{L^2_{t} L^{2}_x}.\eneq Moreover, if
$F\equiv 0$ and $s\in[0,1]$, then \beeq\label{70-est-LocEnerDecay} \big\Vert
\beta u \big\Vert_{L^2_{t}H^s}\lesssim    \Vert
u_0\Vert_{\dot{H}^{s}}+\Vert u_1\Vert_{\dot H^{s-1}_x}.\eneq
\end{lem}
\begin{prf}
By \eqref{70-est-Morawetz2}, noting the support property of the
forcing term, we know that
$$\big\Vert\<x\>^{-1/2-\vep} (\pt, \pa_x ) u \big\Vert_{L^2_{t,x}} +\big\Vert \<x\>^{-3/2-\vep} u \big\Vert_{L^2_{t,x}}
\lesssim    \Vert u_0\Vert_{\dot{H}^{1}}+\Vert u_1\Vert_{L^2_x}+
\Vert F\Vert_{L^2_{t,x} }.$$ Thus
$$\big\Vert \beta u \big\Vert_{L^2_{t}H^1}\les\big\Vert \<x\>^{-1/2-\vep} (u, \pa_x u)
\big\Vert_{L^2_{t,x}}\les \Vert u_0\Vert_{\dot{H}^{1}}+\Vert
u_1\Vert_{L^2_x}+ \Vert F\Vert_{L^2_{t,x} }.$$

In the case $F=0$, let $v=g u$, then we have
$(\pt^2-g\Delta_{\mathfrak{g}}g^{-1})v=0$. Note that $\|g u\|_{\dot
H^1}\les \|u\|_{\dot H^1}$, this gives $\|g u\|_{\dot H^{-1}}\les
\|u\|_{\dot H^{-1}}$ by duality. Thus we have by
\eqref{70-est-Morawetz},
$$\big\Vert \beta u \big\Vert_{L^2_{t,x}} \lesssim \big\Vert \beta v\big\Vert_{L^2_{t,x}}
\les\big\Vert \<x\>^{-1/2-\ep} v \big\Vert_{L^2_{t,x}} \les  \Vert
v(0,\cdot)\Vert_{L^2}+\Vert \pt v(0,\cdot)\Vert_{\dot{H}^{-1}}\les
\Vert u_0\Vert_{L^2}+\Vert u_1\Vert_{\dot{H}^{-1}}.$$ This completes
the proof by interpolation with \eqref{70-est-LocEnerDecay1}.
\end{prf}

Now we are ready to give the proof of the KSS estimates presented in
Theorem \ref{70-thm-KSS}.
In the proof, we shall cut the solution into two parts:
a spatially localized part and the part near spatial infinity.
For the localized part, we can use the local energy decay in Lemma
\ref{70-thm-LocEnerDecay}, while for the part near infinity, we can view
the metric as small perturbation of the Minkowski metric and use
Lemma \ref{70-thm-KSS-SmallPert}.

To this end, we introduce a cutoff function $\phi\in C^{\infty}_c$
which equals $1$ in the unit ball $B_1$ and support in $B_{2}$ and
let $\phi_R(x)=\phi(x/R)$. Let $v=\phi_R u$ and $w=(1-\phi_R)u$ with
$R\gg 1$ to be determined later, then
$$\Box_{\mathfrak{g}} v=[-\Delta_\mathfrak{g},\phi_R]u+\phi_R F:= f + \phi_R F.$$
$$\Box_{\mathfrak{g}} w=-[-\Delta_\mathfrak{g},\phi_R]u+(1-\phi_R) F:= -f + (1-\phi_R) F.$$

By \eqref{70-est-Morawetz2}, we have
\beeq\label{70-est-Est-Commutator}\|f\|_{L^2_{t,x}}\les \|\phi_{2R}
\pa_x u\|_{L^2_{t,x}}+\|\phi_{2R} u\|_{L^2_{t,x}} \les \Vert
u_0\Vert_{\dot{H}_x^{1}}+\Vert u_1\Vert_{L^2_x}+ \|\<x\>^{1/2+\ep} F
\|_{L^2_{T,x}}.\eneq
 Note that
$v$ is compactly supported in $x$, so Lemma
\ref{70-thm-LocEnerDecay} applies, i.e.
\beeq\label{70-est-Est-LocalSlt} \big\Vert v \big\Vert_{L^2_t H^1_x}
+\big\Vert \pt v \big\Vert_{L^2_{t,x}}\les \Vert
u_0\Vert_{\dot{H}^{1}_x}+\Vert u_1\Vert_{L^2_x}+ \|\<x\>^{1/2+\ep}
F\|_{L^2_{T,x}}.\eneq

If we define \beeq\label{70-def-Gtilde}\tilde{g}^{ij}=(1-\phi_{R/2})g^{ij}+\phi_{R/2} \delta^{ij},\eneq then
$$\tilde{\Delta}=  \sum_{ij} {\tilde g^{-2}} \partial_i \tilde g^{ij} \tilde g^2 \partial_j =
\Delta+h^{ij} \partial_i \partial_j +b^i\pa_i+c$$ where $h^{ij}=(1-\phi_{R/2})( g^{ij}-\delta^{ij})$,
$b^i=(1-\phi_{R/4}) r^i_1$ and $c=(1-\phi_{R/4})  r_2$, with the $r_j$ satisfying
the bounds in \eqref{c17}. Thus since, $\tilde{\Delta} w=\Delta_\mathfrak{g} w$, $w$ satisfies
\beeq\label{70-est-Est-EqnForOuter} ( \pt^2 - \Delta - h^{ij} ) w =
\Box_\mathfrak{g} w + (b^i \pa_i +c) w = (1-\phi_R) F - f + (b^i
\pa_i +c)w :=G.\eneq

Note that by \eqref{c1} and choosing $R\gg 1$ large enough, $h^{ij}$ will be small enough so that we can apply
Lemma \ref{70-thm-KSS-SmallPert} with $\ep<\rho/4$,
\begin{eqnarray*}
\textrm{RHS of } \eqref{70-est-KSS-SmallPert} & \les & \Vert w' ( 0 , \cdot ) \Vert^2_{L^2 ( \R^{d} )} +
\int^T_0 \int \left(|w'|+\frac{w}{|x|}\right) \left(|G| +(|h'|+\frac{h}{|x|})|w'|\right)d x dt\\
&\les& \Vert w' ( 0 , \cdot ) \Vert^2_{L^2 ( \R^{d} )} +  \|\<x\>^{-1/2-\ep} \left(|w'|+\frac{w}{|x|}\right)\|_{L^2_{T,x}}
\|\<x\>^{1/2+\ep} G\|_{L^2_{T,x}}\\
&&+  \|\<x\>^{-1/2-\rho/2}\left(|w'| + \frac{w}{|x|} \right)\|_{L^2_{T,x}}^2\\
&\les&\Vert w' ( 0 , \cdot ) \Vert^2_{L^2 ( \R^{d} )} +  \vep\|\<x\>^{-1/2-\ep}
\left(|w'|+\frac{w}{|x|}\right)\|^2_{L^2_{T,x}} + \frac{1}{\vep} \|\<x\>^{1/2+\ep}G\|^2_{L^2_{T,x}}\\
&&+ R^{-\rho/2} \|\<x\>^{-1/2-\rho/4}\left(|w'|+\frac{w}{|x|}\right) \|_{L^2_{T,x}}^2 \\
& \les & \Vert w' ( 0 , \cdot ) \Vert^2_{L^2 ( \R^{d} )} + 2 \vep
\|\<x\>^{-1/2-\ep} \left( |w'| + \frac{w}{|x|} \right)
\|^2_{L^2_{T,x}} + \frac{1}{\vep} \| \<x\>^{1/2+\ep} G
\|^2_{L^2_{T,x}}.
\end{eqnarray*}
Since we know from Lemma \ref{70-est-Morawetz-LowOrder},
$$ \| \<x\>^{-1/2-\ep} \left( |w'| + \frac{|w|}{|x|} \right) \|_{L^2_{T,x}} \les \| \<x\>^{-1/2-\ep} |u'| +
\<x\>^{-3/2-\ep} |u| \|_{L^2_{T,x}} \les \| u'(0) \|_{L^2_x} + \|
\<x\>^{1/2+\ep} F \|_{L^2_{T,x}},$$ and
$$ \| \<x\>^{1/2+\ep} (b^i \pa_i +c) w \|_{L^2_{T,x}} \les \| \<x\>^{-1/2+\ep-\rho} \left( |w'| + \frac{w}{|x|} \right)
\|_{L^2_{T,x}} \les \| \<x\>^{-1/2-\ep}
\left(|w'|+\frac{w}{|x|}\right) \|_{L^2_{T,x}},$$ Thus combining
\eqref{70-est-Est-Commutator},
$$\textrm{RHS of }\eqref{70-est-KSS-SmallPert} \les \Vert u_0 \Vert_{\dot{H}^{1}} + \Vert u_1 \Vert_{L^2_x}
+ \| \<x\>^{1/2+\ep} F \|_{L^2_{T,x}}.$$ Applying Lemma
\ref{70-thm-KSS-SmallPert} and \eqref{70-est-Est-LocalSlt}, we get
finally that
$$
(\ln(2+T))^{-1/2} \big\Vert \< x \>^{-1/2} (|u'| +
\frac{|u|}{\<x\>}) \big\Vert_{L^2_{T,x}}
\les \Vert u_0\Vert_{\dot{H}^{1}}+\Vert u_1\Vert_{L^2_x}+
\|\<x\>^{1/2+\ep} F\|_{L^2_{T,x}},
$$
Therefore, we can use Lemma \ref{70-est-Morawetz-LowOrder} to see that $w$ also satisfies
the bounds in \eqref{70-est-KSS}, which completes the proof of Theorem \ref{70-thm-KSS}.

\section{Higher Order KSS Estimates}\label{70-sec-Kss-high}
In this section, we prove the KSS estimates involving high order derivatives as stated in Theorem \ref{70-thm-KSS-HighOrder}.

We first give the KSS estimate for $P$. Consider the equation
\begin{equation} \label{70-eqn-LW-transf}
\begin{cases}
(\partial_t^2+P)v(t,x)=G(t,x),
\quad (t,x)\in \R_+\times \R^d
\\
u(0,x)=v_0(x), \quad \partial_t v(0,x)=v_1(x), \quad x\in \R^d.
\end{cases}
\end{equation}
Recall that if we let $v=g u$ and $G=g F$, then \beeq\label{70-est-EquiEqn}(\pt^2-\Delta_{\mathfrak{g}})u=F \Leftrightarrow (\pt^2+P)v=G.\eneq
Thus we have also the following KSS estimate
\begin{equation} \label{70-est-KSSv}
A_{\mu}(T)\left(\big\Vert \< x \>^{-\mu} v' \big\Vert_{L^2_{T,x}}+\big\Vert \< x \>^{-\mu-1} v \big\Vert_{L^2_{T,x}} \right) \lesssim    \Vert v' ( 0 , \cdot ) \Vert_{L^2_x} +  \Vert \< x \>^{\frac{1}{2}+\vep} G\Vert_{L^2_{T,x}}
\end{equation} for $\mu\ge 1/2$ and $\ep>0$.

For $\pa_x^\al \Omega^\be v$, we will use induction on $|\be|$ to give the proof. First, for $|\be|=0$, we have $$(\pa_t^2+P)\pa_x^\al v=\pa_x^\al  G+[P,\pa_x^\al]v.$$
Locally, by ellipticity, we have
 $$\sum_{|\al|=2}\|\pa_x^\al f\|_{L^2_{|x|\le R}}\les \sum_{|\be|\le 1} \| \pa^\be_x f\|_{L^2_{|x|\le 2R}}+\|P f\|_{L^2_{|x|\le 2R}},$$
and so,  by an induction argument, and \eqref{a} and \eqref{b} in
Lemma \ref{70-lem-equi-PLapBH},
 \begin{eqnarray*}
\|\pa_x^\al f\|_{L^2_{|x|\le R}} &\les&
\sum_{2 j+1\le |\al|} \|\pa_x P^j f\|_{L^2_{|x|\le 2R}}+\sum_{2 j\le |\al|} \|P^j f\|_{L^2_{|x|\le 2R}}\\
&\les& \| \<x\>^{-2} f\|_{L^2} + \sum_{j\le |\al|-1}
\|\<x\>^{-1}\pa_x P^{j/2}  f\|_{L^2}. \end{eqnarray*} Thus by the
KSS estimates \eqref{70-est-KSSv}, \eqref{b} and \eqref{e},
\begin{eqnarray}
  \sum_{|\al|\le N+1} \| \pa_x^\al v\|_{L^2_{t,|x|\le R}} &\les& \sum_{j\le N} \|\<x\>^{-1} \pa_x P^{j/2} v\|_{L^2_{t,x}} + \| \<x\>^{-2} v\|_{L^2}\nonumber\\
  &\les& \sum_{j\le N}\left( \|\pa_x P^{j/2} v_0\|_{L^2}+\|P^{j/2} v_1\|_{L^2}+\|\<x\>^{1/2+\ep/2} P^{j/2} G \|_{L^2_{t,x}} \right) \nonumber\\
  &\les& \sum_{\be\le N}\left( \|\pa_x^\be v_0\|_{\dot H^1}+\|\pa_x^\be v_1\|_{L^2}+\|\<x\>^{1/2+\ep} \pa_x^\be G \|_{L^2_{t,x}}\right)  .
\end{eqnarray}

Using again \eqref{70-est-KSSv}, we have
\begin{align}
  A_\mu(T) \sum_{|\al|\le N}\big\Vert \< x \>^{-\mu} (|(\pa_x^\al v)'|+\frac{|\pa_x^\al  v|}{\<x\>}) \big\Vert_{L^2_{T,x}}\les\label{70-est-KSS-HiOrToGet}\\
  \sum_{|\al|\le N}\left( \Vert \pa_x^\al  v_0\Vert_{\dot{H}^{1}}+\Vert \pa_x^\al v_1\Vert_{L^2_x}+
\|\<x\>^{1/2+\ep}\pa_x^\al  G\|_{L^2_{T,x}}+\|\<x\>^{1/2+\ep}[P,\pa_x^\al] v\|_{L^2_{T,x}}\right)\nonumber\end{align} for $\mu\ge 1/2$ and $\ep>0$.
Note that $$ |[P,\pa_x^\al]v| \les \<x\>^{-\rho-1} \sum_{1\le |\ga|\le |\al|+1} |\pa_x^{\ga} v| + \<x\>^{-\rho-2} |v|, $$ if choose $\ep>0$ small enough
so that $1/2+\ep-\rho-1<-1/2-3 \ep$ and $R\gg 1$,
$$\sum_{|\al|\le N}\|\<x\>^{1/2+\ep}[P,\pa_x^\al] v\|_{L^2_{T,x}}\les
\vep \sum_{|\al|\le N}\|\<x\>^{-1/2-\ep} (|(\pa_x^\al
v)'|+\frac{|\pa_x^\al v|}{\<x\>})\|+ \sum_{|\al|\le N+1}\|\pa_x^\al
v\|_{L^2_t L^2_{|x|\le R}}.$$ So the first term in the right hand
side of the inequality can be absorbed by
\eqref{70-est-KSS-HiOrToGet} with $\mu=1/2+\ep$. Thus for $\rho>0$,
$\mu\ge 1/2$ and $\ep>0$, \beeq\label{70-est-KSS-HighOrder-x}
A_\mu(T) \sum_{|\al|\le N}\big\Vert \< x \>^{-\mu}
(|(\pa_x^\al v)'|+\frac{|\pa_x^\al v|}{\<x\>})
\big\Vert_{L^2_{T,x}}\les
  \sum_{|\al|\le N} \left( \Vert \pa_x^\al  v_0\Vert_{\dot{H}^{1}}+\Vert \pa_x^\al v_1\Vert_{L^2_x}+
\|\<x\>^{1/2+\ep} \pa_x^\al  G\|_{L^2_{T,x}}\right).\eneq

Now we claim that for $\rho>1$, $\mu\ge 1/2$ and $\ep>0$, we have \begin{align}
  A_\mu(T) \sum_{|\al|+|\be|\le N}\big\Vert \< x \>^{-\mu} (|(\pa_x^\al \Omega^\be v)'|+\frac{|\pa_x^\al \Omega^\be  v|}{\<x\>}) \big\Vert_{L^2_{T,x}}\les\nonumber\\
  \sum_{|\al|+|\be|\le N}\left( \Vert \pa_x^\al \Omega^\be  v_0\Vert_{\dot{H}^{1}}+\Vert \pa_x^\al \Omega^\be v_1\Vert_{L^2_x}+
\|\<x\>^{1/2+\ep}\pa_x^\al \Omega^\be  G\|_{L^2_{T,x}}\right).\label{70-est-KSS-HighOrderInduct}\end{align}
We use induction on $|\be|$ to give the proof. Assume the estimate \eqref{70-est-KSS-HighOrderInduct} is true for $|\be|\le k$, then
$$ (\pa_t^2+P) \Omega v = \Omega G + [P,\Omega] v $$
and we are reduced to estimate
$$\sum_{|\al|+|\be|\le N-1, |\be|\le k}\|\<x\>^{1/2+\ep}\pa_x^\al \Omega^\be  [P,\Omega] v\|_{L^2_{T,x}}$$
by \eqref{70-est-KSS-HighOrderInduct}. Note that
$$[P, \Omega]v =\sum_{|\ga|\le 2}  r_{|\ga|-2}\pa_x^{\ga} v ,$$ where $r_j$ are functions satisfying
\eqref{c17}. Thus if $\rho>1$ and $\ep>0$ small enough such that $\rho>1+2\ep$, \begin{eqnarray*}
  &&\sum_{|\al|+|\be|\le N-1, |\be|\le k}\|\<x\>^{1/2+\ep}\pa_x^\al \Omega^\be  [P,\Omega] v\|_{L^2_{T,x}}\\
  &\les & \sum_{|\al|+|\be|\le N, |\be|\le k}\|\<x\>^{1/2+\ep-\rho} (|(\pa_x^\al \Omega^\be  v)'|+\frac{|v|}{\<x\>})\|_{L^2_{T,x}}\\
  &\les & \sum_{|\al|+|\be|\le N, |\be|\le k}\|\<x\>^{-1/2-\ep} (|(\pa_x^\al \Omega^\be  v)'|+\frac{|v|}{\<x\>})\|_{L^2_{T,x}}\les \eqref{70-est-KSS-HighOrderInduct}.
\end{eqnarray*}

Since $\pt$ commutate with $\pt^2+P$, we can
conclude that for $\rho>1$, $\mu\ge 1/2$ and $\ep>0$, we have \begin{align}
  A_\mu(T) \sum_{j+|\al|+|\be|\le N}\big\Vert \< x \>^{-\mu} (|(\pt^j \pa_x^\al \Omega^\be v)'|+\frac{|\pt^j \pa_x^\al \Omega^\be  v|}{\<x\>}) \big\Vert_{L^2_{T,x}}\les\nonumber\\
  \sum_{j+|\al|+|\be|\le N}\left( \Vert (\pt^j \pa_x^\al \Omega^\be v)'(0,\cdot)\Vert_{L^{2}}+
\|\<x\>^{1/2+\ep}\pt^j \pa_x^\al \Omega^\be  G\|_{L^2_{T,x}}\right).\label{70-est-KSS-HighOrderAll}\end{align}
Combining the energy estimates, we get the following estimates
\begin{align*}
\sup_{0 \leq t \leq T} \sum_{1 \leq k + j \leq N+1} \big\Vert \partial_{t}^{k} P^{j/2} v (t & , \cdot ) \big\Vert_{L^2_x} + \sum_{\vert \alpha \vert \leq N} A_{\mu}(T) \big\Vert \<x\>^{-\mu} (|(Z^{\alpha} v)'|+  \frac{|Z^{\alpha} v|}{\<x\>}) \big\Vert_{L^2_T L^2_x} \nonumber \\
& \lesssim \sum_{\vert \alpha \vert \leq N}\bigg( \big\Vert (Z^{\alpha} v)' (0, \cdot ) \big\Vert_{L^2_x} +  \int_0^T \big\Vert Z^{\alpha} G(s, \cdot ) \big\Vert_{L^2_x} d s \bigg).
\end{align*} for the solution to the equation \eqref{70-eqn-LW-transf}. For the $\big\Vert (Z^{\alpha} v)' (0, \cdot ) \big\Vert_{L^2_x}$ part, if it has $\pt^j$ component with $j\ge 2$, we can use the equation to reduce it to the case $\pt^{j-2}$, with an additional term $$\sum_{\vert \alpha \vert \leq N-1} \|Z^\al G(0,\cdot)\|_{L^2_x} \les
\sum_{\vert \alpha \vert \leq N-1} \|Z^\al G(t,\cdot)\|_{L^2_x W^{1,1}_t}\les \sum_{\vert \alpha \vert \leq N} \|Z^\al G(t,\cdot)\|_{L^1_T L^2_x}.$$ This means that
\begin{align}\label{70-est-KSS-HighOrderv}
\sup_{0 \leq t \leq T} \sum_{1 \leq k + j \leq N+1} \big\Vert \partial_{t}^{k} P^{j/2} v (t & , \cdot ) \big\Vert_{L^2_x} + \sum_{\vert \alpha \vert \leq N} A_{\mu}(T) \big\Vert \<x\>^{-\mu} (|(Z^{\alpha} v)'|+  \frac{|Z^{\alpha} v|}{\<x\>}) \big\Vert_{L^2_T L^2_x} \nonumber \\
& \lesssim \sum_{\vert \alpha \vert \leq N} \big\Vert (Y^{\alpha} v)' (0, \cdot ) \big\Vert_{L^2_x} + \sum_{\vert \alpha \vert \leq N} \int_0^T \big\Vert Z^{\alpha} G(s, \cdot ) \big\Vert_{L^2_x} d s .
\end{align}

Turning back to the equation \eqref{70-eqn-LW}, let $v=g u$ and $G=g F$, then $(\pt^2+P^2)v=G$ with $v_0=g u_0$, $v_1=g u_1$. Note that $$\pa_x u=g^{-1}\pa_x g u-g^{-1}(\pa_x g) u=g^{-1}\pa_x v-g^{-2}(\pa_x g) v,$$
we have $$\sum_{|\alpha| \leq N} \big\Vert \<x\>^{-\mu} (|(Z^{\alpha} u)'|+  \frac{|Z^{\alpha} u|}{\<x\>}) \big\Vert_{L^2_{T,x}}\les \sum_{|\alpha| \leq N} \big\Vert \<x\>^{-\mu} (|(Z^{\alpha} v)'|+  \frac{|Z^{\alpha} v|}{\<x\>}) \big\Vert_{L^2_{T,x}}.$$
This concludes the proof of Theorem \ref{70-thm-KSS-HighOrder}.

\section{Almost Global Existence}\label{70-sec-AlmostGlob}

In this section we shall prove one of our main existence theorems, Theorem \ref{70-thm-DSLW-AGlob}. The proof will be similar to that of Keel, Smith and Sogge for the Minkowski case (see \cite{KeSmSo02}). We start with the now standard Sobolev estimate (see \cite{Kl85_01}).

\begin{lem}\sl   \label{LN1}
Suppose that $h \in C^{\infty} ( \R^3 )$. Then, for $R>1$,
\begin{equation} \label{N1}
\Vert h \Vert_{L^{\infty} (R/2 \leq \vert x \vert \leq R)} \lesssim R^{-1} \sum_{\vert \alpha \vert \leq 2} \Vert Y^{\alpha} h \Vert_{L^2(R/4 \leq \vert x \vert \leq 2R)}.
\end{equation}
\end{lem}

We now define the bilinear form $\widetilde{Q}$ by $\widetilde{Q}(u',u')=Q(u')$. The following estimate for the nonlinear part will be crucial.

\begin{lem}\sl  \label{LN2}
We have
\begin{equation*}
\sum_{|\al|\le 4}\big\Vert Z^{\al} \widetilde{Q} (u',v')
\big\Vert_{L^2 ( \R^{3} )}^{2} \lesssim \bigg( \sum_{\vert \alpha
\vert \leq 4} \big\Vert \< x \>^{- 1/2} Z^{\alpha} u'
\big\Vert^2_{L^2 ( \R^{3} )} \bigg) \bigg( \sum_{\vert \alpha \vert
\leq 4} \big\Vert \< x \>^{- 1/2} Z^{\alpha} v' \big\Vert^2_{L^2 (
\R^{3} )} \bigg) .
\end{equation*}
\end{lem}

\begin{prf}
We clearly have the pointwise bound:
\begin{align*}
\big\vert Z^{\beta} \widetilde{Q} (u',v') (s,x) \big\vert \lesssim & \bigg( \sum_{\vert\alpha\vert\leq 4} \big\vert Z^{\alpha}u'(s,x) \big\vert \bigg) \bigg( \sum_{\vert \alpha \vert \leq  2 }
 \big\vert Z^{\alpha} v' (s,x) \big\vert \bigg)   \\
&+ \bigg( \sum_{\vert\alpha\vert\leq 4} \big\vert Z^{\alpha} v'(s,x)
\big\vert \bigg) \bigg( \sum_{\vert \alpha \vert \leq 2} \big\vert
Z^{\alpha} u' (s,x) \big\vert \bigg) .
\end{align*}
We need only to estimate the first term. Using Lemma \ref{LN1} for a given $R=2^j$, $j\geq 0$, we get
\begin{align*}
\big\Vert Z^{\beta} & \widetilde{Q} (u',v') \big\Vert^{2}_{L^2 ( \{ \vert x \vert \in [ 2^j , 2^{j+1} ] \} )}  \\
&\lesssim 2^{-2j } \sum_{\vert \alpha \vert \leq 4} \big\Vert Z^{\alpha} u' \big\Vert^{2}_{L^2 ( \{ \vert x \vert \in [ 2^j , 2^{j+1} ] \} )} \sum_{\vert \alpha \vert \leq 4} \big\Vert Z^{\alpha} v' \big\Vert^{2}_{L^2 ( \{ \vert x \vert \in [ 2^{j-1} , 2^{j+2} ] \} )}    \\
&\lesssim \sum_{\vert \alpha \vert \leq 4} \big\Vert \< x \>^{- 1/2} Z^{\alpha} u' \big\Vert^{2}_{L^2 ( \{ \vert x \vert \in [ 2^{j} , 2^{j+1} ] \} )} \sum_{\vert \alpha \vert \leq 4}
 \big\Vert \< x \>^{- 1/2} Z^{\alpha} v' \big\Vert^{2}_{L^2 ( \{ \vert x \vert \in [ 2^{j-1} , 2^{j+2} ] \} )}   \\
&\lesssim \sum_{\vert \alpha \vert \leq 4} \big\Vert \< x \>^{- 1/2}
Z^{\alpha} u' \big\Vert^{2}_{L^2 ( \{ \vert x \vert \in [ 2^{j} ,
2^{j+1} ] \} )}  \sum_{\vert \alpha \vert \leq 4} \big\Vert \< x
\>^{- 1/2} Z^{\alpha} v' \big\Vert^{2}_{L^2 ( \R^{3} )} .
\end{align*}
We also have the bound
\begin{eqnarray*}
\big\Vert Z^{\beta} \widetilde{Q} (u',v') \big\Vert_{L^2 ( \{ \vert x \vert < 1 \} )}^{2} \lesssim \sum_{\vert \alpha \vert \leq L} \big\Vert Z^{\alpha} u' \big\Vert^2_{L^2 ( \{ \vert x \vert <2 \} )} \sum_{\vert \alpha \vert \leq L} \big\Vert Z^{\alpha} v' \big\Vert^2_{L^2 ( \{ \vert x \vert <2 \} )} .
\end{eqnarray*}
Summing over $j$ gives the lemma.
\end{prf}

{\bf\noindent Proof of Theorem \ref{70-thm-DSLW-AGlob}.}
We follow \cite{KeSmSo02}. Let $u_{-1}=0$. We define $u_k$, $k \in \N$ inductively by letting $u_k$ to solve
\begin{equation}   \label{N3}
\begin{cases}
\Box_{\mathfrak{g}} u_k=Q ( u_{k-1}'),
\\
u(0,x)=u_0(x), \quad \partial_t u(0,x)=u_1(x).
\end{cases}
\end{equation}
For $T > 0$, we denote
\begin{equation*}
M_k (T) = \sup_{0 \leq t \leq T} \sum_{1 \leq i + j \leq 5} \big\Vert \partial_{t}^{i} P^{j/2} g u_k \big\Vert_{L^{2} ( \R^{3} )} + \sum_{\vert \alpha \vert \leq 4} (\ln (2+T))^{-1/2} \big\Vert \< x \>^{- 1/2} Z^{\alpha} u_k' \big\Vert_{L^2( [0 , T] \times \R^{3} )} .
\end{equation*}
Using Theorem \ref{70-thm-KSS-HighOrder}, we see that there exists a constant $C_0$ such that
\begin{equation*}
M_0 (T) \leq C_0 \delta ,
\end{equation*}
for any $T$. We claim that, for $k \geq 1$, we have
\begin{equation} \label{c10}
M_k ( T_{\delta} ) \leq 2 C_0 \delta ,
\end{equation}
for $\delta$ sufficiently small and $T_{\delta}$ appropriately chosen later. We will prove this inductively. Assume that the bound holds for $k-1$. By Theorem \ref{70-thm-KSS-HighOrder}, we have, for $\delta$ small enough,
\begin{align*}
M_k ( T_{\delta} ) &\leq C_0 \delta + C \sum_{\vert \alpha \vert \leq 4} \int_0^{T_{\delta}} \big\Vert Z^{\alpha} Q ( u_{k-1} ') (s, \cdot ) \big\Vert_{L^2 ( \R^d )} d s   \\
& \leq C_0 \delta + C \sum_{\vert \alpha \vert \leq 4} \int_0^{T_{\delta}} \big\Vert \< x \>^{- 1/2} Z^{\alpha} u_{k-1}' \big\Vert^2_{L^2 ( \R^d )} d s   \\
& \leq C_0 \delta + C \ln (2+T_\de) M_{k-1}^2 ( T_{\delta} )   \\
& \leq C_0 \delta + C \ln (2+T_\de)( 2 C_0 \delta)^2 ,
\end{align*}
where we have also used Lemma \ref{LN2} and the induction hypothesis. Then, to prove \eqref{c10}, it is enough to have
\begin{equation} \label{N4}
C_0 \delta + C \ln (2+T_\de) ( 2 C_0 \delta )^2 \leq 2 C_0 \delta \Longleftrightarrow 4 C C_0 \ln (2+T_\de) \delta \leq 1 .
\end{equation}
Therefore, we can set $T_{\delta} = \exp(c/\de)$ and $c$ small enough.

To show that the sequence $u_k$ converges, we estimate the quantity
\begin{align*}
A_k (T) =  & \sup_{0 \leq t \leq T} \sum_{1 \leq i + j \leq 5} \big\Vert \partial_{t}^{i} P^{j/2} g ( u_k - u_{k-1} ) \big\Vert_{L^2 ( \R^d )}  \\
& + \sum_{\vert \alpha \vert \leq 4} (\ln (2+T_\de))^{-1/2} \big\Vert \< x \>^{- 1/2} Z^{\alpha} ( u'_k - u'_{k-1} ) \big\Vert_{L^2( [0,T] \times \R^{d} )}.
\end{align*}
It is clearly sufficient to show
\begin{equation} \label{N5}
A_k (T) \leq \frac{1}{2} A_{k-1} (T) .
\end{equation}
Using Lemma \ref{LN2} and repeating the above arguments, we obtain
\begin{align*}
A_k ( T_{\delta} ) & \leq  \widetilde{C} \sum_{\vert \alpha \vert \leq 4} \int_0^{T_{\delta}}  \big\Vert \< x \>^{- 1/2} Z^{\alpha} ( u'_{k-1} - u'_{k-2} ) \big\Vert_{L^2 ( \R^d )}  \\
& \times \sum_{\vert \alpha \vert \leq 4} \Big( \big\Vert \< x \>^{- 1/2} Z^{\alpha} u'_{k-1} \big\Vert_{L^2 ( \R^d )} + \big\Vert \< x \>^{- 1/2} Z^{\alpha} u'_{k-2} \Vert_{L^2 ( \R^d )} \Big) d s \\
& \leq \widetilde{C} \ln(2+ T_{\delta} ) ( M_{k-1} ( T_{\delta} ) + M_{k-2} ( T_{\delta} ) ) A_{k-1} (T_{\delta} ) .
\end{align*}
Using \eqref{c10}, the above inequality leads to \eqref{N5} if $\delta$ is small enough.
{\hfill  {\vrule height6pt width6pt depth0pt}\medskip}

\section{Global Strichartz Estimates}\label{70-sec-Strichartz}
In this section we shall prove Theorem \ref{70-thm-Strichartz-Glob}, the global Strichartz estimates for $\Box_{\mathfrak{g}}$. The proof will follow closely the arguments of \cite{SmSo00}, \cite{Burq} and \cite{HMSSZ}, which dealt with compact perturbations of the Minkowski metric or obstacle case.

Recall that we have equivalence \eqref{70-est-EquiEqn} between the wave equations involving $P$ and $-\Delta_\mathfrak{g}$. We need only to give the proof for the case of following equation,
\begin{equation} \label{70-eqn-LW-Stri}
\begin{cases}
(\partial_t^2+P)u(t,x)=F(t,x),
\quad (t,x)\in \R_+\times \R^d
\\
u(0,x)=u_0(x), \quad \partial_t u(0,x)=u_1(x), \quad x\in \R^d.
\end{cases}
\end{equation}

As in \cite{HMSSZ}, we will need the following four ingredients: local energy decay for $\pt^2+P$, local Strichartz estimates for $\pt^2+P$, global Strichartz estimates for $\pt^2+\tilde{P}$ with $\tilde{\mathfrak{g}}=\mathfrak{g}$ in $|x|\ge R$ for some $R>0$, and an estimate for $\pt^2+\tilde{P}$ like Proposition 2.1 in \cite{HMSSZ}.

Let us give first the estimates for $\pt^2+\tilde{P}$, where $\tilde{\mathfrak{g}}$ is defined by \eqref{70-def-Gtilde}, i.e., $\tilde{g}^{ij}=(1-\phi_{R/2})g^{ij}+\phi_{R/2} \delta^{ij}$.
Then $$\tilde{P}= - \sum_{ij} \frac{1}{\tilde g} \partial_i \tilde g^{ij} \tilde g^2 \partial_j \frac{1}{\tilde g}=-\pa_i \tilde g^{ij} \pa_j +b^i\pa_i+c$$ where $b^i=-2 \tilde{g}^{-1} \pa_j(\tilde g \tilde g^{ij})-\tilde g \tilde g^{ij} \pa_j \tilde g^{-1}+\pa_j \tilde g^{ij}$ and $c=\tilde g^{-1} \pa_i (\tilde g^{ij}\pa_j \tilde g)$.

Consider the wave equation
\begin{equation} \label{70-eqn-LW2}
\begin{cases}
(\pt^2+\tilde{P})u=0,
\quad (t,x)\in \R_+\times \R^d
\\
u(0,x)=u_0(x), \quad \partial_t u(0,x)=u_1(x), \quad x\in \R^d.
\end{cases}
\end{equation}

We denote
\[
A_j = \R \times \{ 2^j \leq |x| \leq 2^{j+1}\}, \qquad A_{<j} =  \R \times \{|x| \leq 2^{j}\}.
\]
By hypotheses \eqref{c1}, we can see that if $R\gg 1$ so that $R^{-\rho}\le \ep$,
\begin{equation}
  \label{coeff}
  \sum_{j \in \Z} \sup_{A_j} |x|^2 |\nabla^2 \tilde{g}^{ij}(x)| + |x| |\nabla \tilde{g}^{ij}(x)|
  + |\tilde{g}^{ij}(x)-\delta^{ij}| \les \sum_{j \gtrsim \ln(R)} \sup_{A_j} |x|^{-\rho} \les R^{-\rho}\le \ep
\end{equation}
and, for the lower order terms,
\begin{equation}
  \label{coeffb}
  \sum_{j \in \Z} \sup_{A_j}  |x|^2 |\nabla b^i(x)|
  + |x| |b^i(x)| \les \sum_{j \gtrsim \ln(R)} \sup_{A_j} |x|^{-\rho} \le \ep,
\end{equation}
\begin{equation}
  \label{coeffcc}
  \sum_{j \in \Z} \sup_{A_j}  |x|^4 |c(t,x)|^2  \les \sum_{j \gtrsim \ln(R)} \sup_{A_j} |x|^{-2\rho} \le \ep.
\end{equation}

Thus we can apply Theorem 6 in Metcalfe-Tataru \cite{MeTa07p} to get the following
\begin{prop}[Global Strichartz Estimates for Small Perturbation]\label{70-prop-StriSmallPet}
Let $d\ge 3$, $s\in [0,1]$ ($s\in (0,1)$ if $d=3$). The solution of \eqref{70-eqn-LW2} satisfies
\begin{equation}\label{70-est-StriSmallPet}
\big\Vert u  \big\Vert_{L^q_t L^r_x}\lesssim \|u_0\|_{\dot
H^s}+\|u_1\|_{\dot H^{s-1}},
\end{equation} for any admissible $(s,q,r)$ with $q>2$, $r<\infty$. In other words, we have
\begin{equation}\label{70-est-StriSmallPet2}
\big\Vert e^{i t \tilde{P}^{1/2}} f  \big\Vert_{L^q_t L^r_x}\lesssim \|f\|_{\dot H^s},
\end{equation}
\end{prop}

As a corollary, we can get an estimate for $\pt^2+\tilde{P}$ like Proposition 2.1 in \cite{HMSSZ}.
\begin{prop}\label{70-prop-HMSSZ}
Let $d\ge 3$, $s\in [0,1]$ ($s\in (0,1)$ if $d=3$). The solution $u$ of $(\pt^2+\tilde{P}) u= \beta F=F$ satisfies
\begin{equation}\label{70-est-StriSmallPet}
\big\Vert u  \big\Vert_{L^q_t L^r_x}\lesssim \|u_0\|_{\dot
H^s}+\|u_1\|_{\dot H^{s-1}}+\|F\|_{L^2_t H^{s-1}},
\end{equation} for any admissible $(s,q,r)$ with $q>2$, $r<\infty$.
\end{prop}
\begin{prf}
Applying \eqref{70-est-LocEnerDecay} in this setting, we know that for any $s\in [0,1]$,
$$\|\beta e^{i t \tilde{P}^{1/2}} f\|_{L^2_t H^s}\les \|f\|_{\dot H^s}.$$
By duality, this is equivalent to say that for $s\in [-1,0]$,
$$\|\int  e^{-i t \tilde{P}^{1/2}}\beta F(t,\cdot) dt\|_{\dot H^s}\les \|f\|_{L^2_t H^s}.$$

Combining \eqref{70-est-StriSmallPet2} and Lemma
\ref{70-lem-equi-PLap} \eqref{d},
\begin{eqnarray*}\|\int e^{i (t-s) \tilde{P}^{1/2}} \tilde{P}^{-1/2} \beta F(s,\cdot) d s\|_{L^q_t L^r_x}&\les& \|\int e^{-i s \tilde{P}^{1/2}} \tilde{P}^{-1/2} \beta F(s,\cdot) d s\|_{\dot H^s}\\
&\les&\|\int e^{-i s \tilde{P}^{1/2}} \beta F(s,\cdot) d s\|_{\dot H^{s-1}}\\
&\les&\|F\|_{L^2_t H^{s-1}}.\end{eqnarray*}
Since $q>2$, we can apply the Christ-Kiselev lemma \cite{ChKi01} to conclude the proof.
\end{prf}

Now we give the local Strichartz estimates for the $\pt^2+P$. Consider the wave equation
\begin{equation} \label{70-eqn-LW1}
\begin{cases}
(\pt^2+P) u(t,x)=0,
\quad (t,x)\in \R_+\times \R^d
\\
u(0,x)=u_0(x), \quad \partial_t u(0,x)=u_1(x), \quad x\in \R^d.
\end{cases}
\end{equation}
 The local Strichartz estimates for the variable coefficient wave equations was studied extensively, see e.g. Kapitanski \cite{Ka91}, Mockenhaupt-Seeger-Sogge \cite{MoSeSo93}, Smith \cite{Sm98}, Bahouri-Chemin \cite{BaCh99-1} \cite{BaCh99-2}, Tataru \cite{Ta00} \cite{Ta01} \cite{Ta02}. In particular, we have
\begin{lem}[Theroem 1.1 in Tataru \cite{Ta02}]\label{70-lem-Strichartz}
If $(\pt,\pa_x)^\al a^{ij}(t,x)\in L^1_{t\in [0,1]} L^\infty_x$ for
any $|\al|=2$ and $(\pt^2- \pa_i a^{ij}(t,x) \pa_j)u=F$, then
$$\|D^{1-s} u\|_{L^q_{t\in [0,1]} L^r_x}\les \|u'\|_{L^\infty_{t\in [0,1]} L^2_x}+\|F\|_{L^1_{t\in [0,1]} L^2_x}$$
for any admissible $(s,q,r)$ with $r<\infty$.
\end{lem}

Note that we can write
$$\pt^2+P=\pt^2-\pa_i g^{ij} \pa_j+b^i\pa_i+c$$
with $b^i=\mathcal{O}(\<x\>^{-1-\rho})$ and $c=\mathcal{O}(\<x\>^{-2-\rho})$. Moreover, for $a^{ij}=g^{ij}$ and $|\al|=2$, we have $\pa_x^\al a^{ij}(x)=\mathcal{O}(\<x\>^{-2-\rho})\in L^1_{t\in [0,1]} L^\infty_x$.
Then for the equation \eqref{70-eqn-LW1}, we have
\begin{eqnarray}
  \|D^{1-s} u\|_{L^q_{t\in [0,1]} L^r_x}&\les& \|u'\|_{L^\infty_{t\in [0,1]} L^2_x} + \|b^i\pa_i u+c u\|_{L^1_t L^2_x}\nonumber\\
  &\les& \|u'\|_{L^\infty_{t\in [0,1]} L^2_x}\nonumber\\
  &\les& \|u_0\|_{\dot H^1_x}+\|u_1\|_{L^2_x}\label{70-est-StrichartzLoc1}
\end{eqnarray}
by energy estimates and Lemma \ref{70-lem-Strichartz}. Moreover, we can prove the following

\begin{prop}[Local Strichartz Estimates]\label{70-thm-StrichartzLoc}
Let $s\in [0,1]$, $(s,q,r)$ admissible with $r<\infty$, and $u$ be the solution to the equation \eqref{70-eqn-LW1}, then
\beeq\label{70-est-StrichartzLoc}\|u\|_{L^q_{t\in [0,1]} L^r_x}\les \|u_0\|_{\dot H^s}+\|u_1\|_{\dot H^{s-1}}.\eneq
\end{prop}
\begin{prf}
Since $(\pt^2+P) u=0$, we have $(\pt^2+P) P^{\frac{s-1}{2}}u=0$.
Thus by \eqref{g} and \eqref{d} in Lemma \ref{70-lem-equi-PLap} and
\eqref{70-est-StrichartzLoc1},
\begin{eqnarray*}
  \|u\|_{L^q_{t\in [0,1]} L^r_x}&\les&  \|P^{\frac{1-s}{2}} P^{\frac{s-1}{2}}u\|_{L^q_{t\in [0,1]} L^r_x}\\
  &\les& \|D^{1-s} P^{\frac{s-1}{2}}u\|_{L^q_{t\in [0,1]} L^r_x}\\
  &\les& \|P^{\frac{s-1}{2}}u_0\|_{\dot H^1}+\|P^{\frac{s-1}{2}}u_1\|_{L^2}\\
  &\les& \|u_0\|_{\dot H^s}+\|u_1\|_{\dot H^{s-1}},
\end{eqnarray*}
if $(s,q,r)$ is admissible, $s\in[0,1]$ and $r<\frac{d}{1-s}$, i.e., $s>1-\frac{d}{r}$.
Since $(s,q,r)$ is admissible, then if $d=3$, we must have $q>2$,
$$s=\frac{d}{2}-\frac{d}{r}-\frac{1}{q}>1-\frac{d}{r}.$$
Similarly, if $d\ge 4$, we have $q\ge 2$, $$s=\frac{d}{2}-\frac{d}{r}-\frac{1}{q}\ge \frac{d-1}{2}-\frac{d}{r}>1-\frac{d}{r}.$$
This means that we have the local Strichartz estimates for all the cases where $(s,q,r)$ is admissible and $s\in[0,1]$.
\end{prf}

Now we are ready to give the proof of the global Strichartz estimates for $\pt^2+P$.
\begin{prop}\label{70-thm-Strichartz2nd}
For the linear equation \eqref{70-eqn-LW-Stri}, assume that $$u_0(x)=u_1(x)=F(t,x)=0\quad \textrm{ when }|x|>2 R,$$ then for any $s\in[0,1]$ and admissible $(s,q,r)$ with $q>2$ and $r<\infty$, we have
\beeq\label{70-est-Strichartz2nd} \big\Vert u \big\Vert_{L^q_{t} L^r_x}\lesssim    \Vert u_0\Vert_{\dot{H}^{s}}+\Vert u_1\Vert_{\dot H^{s-1}_x}+  \Vert F \Vert_{L^2_{t} \dot H^{s-1}_x}.\eneq
\end{prop}
\begin{prf}
One of the key ingredients in the proof is the
following variant of Lemma \ref{70-thm-LocEnerDecay}, which holds for all $s\in [0,1]$,
\beeq\label{70-est-est-l2linfty}\|\beta u\|_{L^2_t \dot H^s\cap L^\infty_t \dot H^s}
+\|\beta u_t\|_{L^2_t \dot H^{s-1}\cap L^\infty_t \dot H^{s-1}}
\les \Vert u_0\Vert_{\dot{H}^{s}}+\Vert u_1\Vert_{\dot H^{s-1}_x}+  \Vert F\Vert_{L^2_{t} \dot H^{s-1}_x}.\eneq
The $L^2_t$ estimate with $s=1$ comes from \eqref{70-est-LocEnerDecay1}, then the estimates follow from Duhamel's formula and duality.

To prove \eqref{70-est-Strichartz2nd}, let us argue as before. Let $v=\phi_{3R} u$ and $w=(1-\phi_{3R})u$ with $R\gg 1$ as in the definition of $\tilde{g}$.
Then $w$ solves the wave equation for $\pt^2+\tilde{P}$,
\begin{equation*}
\begin{cases}
(\pt^2+\tilde{P}) w=(\pt^2+{P}) w=[\phi_{3R},P] u
\\
w|_{t=0}=\partial_tw|_{t=0}=0.
\end{cases}
\end{equation*}
An application of Proposition \ref{70-prop-HMSSZ} shows that $\|w\|_{L^q_t L^r_x}$ is
dominated by
$\|\be u\|_{L^2_t \dot H^s}$ if $\be$ equals one on the support of $\phi_{3R}$.
Therefore, by \eqref{70-est-est-l2linfty},
$\|w\|_{L^q_t L^r_x}$ is dominated by the right side of \eqref{70-est-Strichartz2nd}.

As a result, we are left with showing that if $v=\phi_{3R} u$ then
\begin{equation}\label{15}
\|v\|_{L^q_t L^r_x}\lesssim
\|u_0\|_{\dot H^s}+\|g\|_{ \dot H^{s-1}}+\|F\|_{L^2_t \dot H^{s-1}}.
\end{equation}
To do this, fix
$\varphi \in C^\infty_c((-1,1))$ satisfying
$\sum_{j=-\infty}^\infty \varphi(t-j)=1$.  For a given
$j\in {\mathbb N}$, let $v_j=\varphi(t-j)v$.  Then $v_j$ solves
\begin{equation*}
\begin{cases}
(\pt^2+P) v_j = \varphi(t-j)[P ,\phi_{3R}]u -
[\pt^2,\varphi(t-j)]\phi_{3R} u +\varphi(t-j) F
\\
v_j(0,\cdot)=\partial_t v_j(0,\cdot)=0,
\end{cases}
\end{equation*}
while $v_0=v-\sum_{j=1}^\infty v_j$ solves
\begin{equation*}
\begin{cases}
(\pt^2+P) v_0=\tilde{\varphi} [P,\phi_{3R}]u - [\partial_t^2,\tilde{
\varphi}]\phi_{3R} u + \tilde {\varphi }F
\\
v_0|_{t=0}=u_0,\, \, \partial_t v_0|_{t=0}=u_1,
\end{cases}
\end{equation*}
if $\tilde \varphi =1-\sum_{j=1}^\infty \varphi(t-j)$ if $t\ge 0$ and $0$
otherwise.

If we then let $F_j=(\pt^2+P) v_j$ be the forcing term for $v_j$,
$j=0,1,2,\dots$, then, by \eqref{70-est-est-l2linfty},
 we have that
$$\sum_{j=0}^\infty \|F_j\|^2_{L^2_t \dot H^{s-1}}
\lesssim \|u_0\|^2_{\dot H^{s}}+\|u_1\|^2_{\dot H^{s-1}}+
\|F\|_{L^2_t \dot H^{s-1}}^2.$$
By the local Strichartz estimates \eqref{70-est-StrichartzLoc} and Duhamel's formula, we get for
$j=1,2,\dots$
$$\|v_j\|_{L^q_t L^r_x}\lesssim
\int_0^\infty \|F_j(s, \cdot)\|_{ \dot H^{s-1}}\, ds \lesssim
\|F_j\|_{L^2_t \dot H^{s-1}},$$ using Schwarz's inequality and the
support properties of the $F_j$ in the last step.  Similarly,
$$
\|v_0\|_{L^q_t L^r_x}\lesssim \|u_0\|_{\dot
H^s}+\|u_1\|_{H^{s-1}}+\|F_0\|_{L^2_t \dot H^{s-1}}\,.
$$
Since $q, r\ge 2$, we have
$$
\|v\|^2_{L^q_t L^r_x}\lesssim
\sum_{j=0}^\infty \|v_j\|^2_{L^q_t L^r_x}
$$
and so we get
$$
\|v\|^2_{L^q_t L^r_x}\lesssim
\|u_0\|^2_{\dot H^s}+\|g\|^2_{\dot H^{s-1}}+\|F\|_{L^2_t \dot H^{s-1}}^2
$$
as desired, which finishes the proof.
\end{prf}

\noindent{\bf End of Proof of Theorem \ref{70-thm-Strichartz-Glob}:}
Recall that we are assuming that $(\pt^2+P) u=0$.
By Proposition \ref{70-thm-Strichartz2nd}, we may also assume that the initial data
for $u$ vanishes when $|x|<3R/2$.  We then fix $\beta\in C^\infty_c(\R^n)$
satisfying $\beta(x)=1$, $|x|\le R$ and $\beta(x)=0$, $|x|>3R/2$ and write
$$
u=\tilde{u}-v = (1-\beta)\tilde{u} + ( \beta \tilde{u}-v)\,,
$$
where $\tilde{u}$ solves the Cauchy problem for $(\pt^2+\tilde{P}) \tilde{u}=0$
with initial data $(u_0,u_1)$.
By the global Strichartz estimate \eqref{70-est-StriSmallPet},
we can restrict our attention to
$w= \beta \tilde{u}-v$.
But
$$
(\pt^2+P) w = [P, \beta]\tilde{u}\equiv G
$$
is supported in $R<|x|<2R$, and satisfies
\begin{equation}\label{Fdecay}
\int_{0}^\infty \|G(t,\cdot)\|^2_{\dot H^{s-1}}\, dt
\lesssim \|u_0\|_{\dot H^s}^2+\|u_1\|_{\dot H^{s-1}}^2
\end{equation}
by Lemma \ref{70-thm-LocEnerDecay}.
Note also that $w$ has vanishing initial data.
Therefore, since Proposition \ref{70-thm-Strichartz2nd} tells us that
$\|w\|_{L^q_t L^r_x}^2$
is dominated by the left side of \eqref{Fdecay},
the proof is complete.{\hfill  {\vrule height6pt width6pt depth0pt}\medskip}

\section{Weighted Strichartz Estimates}\label{70-sec-W-Strichartz}

Recall Lemma \ref{70-thm-Morawetz} and Remark \ref{70-rem-Morawetz}, we have the following estimates
$$\|\<x\>^{-1/2-\ep}e^{i t P^{1/2}} f\|_{L^2_t L^2_x}\les \|f\|_{L^2_x},$$
$$\|\<x\>^{-3/2-\ep}e^{i t P^{1/2}} f\|_{L^2_t L^2_{x}}\les \|f\|_{\dot H^1_x},$$
By interpolation, we get
\beeq\label{70-est-WSE-L2}\||x|^{-\al-\ep}e^{i t P^{1/2}} f\|_{L^2_t L^2_{|x|\ge R}}\les
\|\<x\>^{-\al-\ep}e^{i t P^{1/2}} f\|_{L^2_t L^2_x}\les
\|f\|_{\dot H^{\al-1/2}_x}\eneq
for any $\al\in [1/2,3/2]$, $\ep>0$  and $d\ge 3$.

Recall that Fang and Wang obtained the following Sobolev
inequalities with angular regularity (Corollary 1.2 in \cite{FaWa})
\beeq\label{70-est-Sobolev-angular} \||x|^{\frac{d}{2}-\al}
f(x)\|_{L^\infty_{|x|} L^{2+\eta}_\omega} \les
\||x|^{\frac{d}{2}-\al} f(x)\|_{L^\infty_{|x|}
H^{\al-\frac{1}{2}}_\omega}
  \les \|
  f\|_{\dot H^{\al}_x} \eneq
for $\al\in (1/2,d/2)$ and some $\eta>0$. Then by \eqref{d} in Lemma
\ref{70-lem-equi-PLap}, we have \beeq\label{70-est-WSE-Infty}
\||x|^{\frac{d}{2}-\al} e^{i t P^{1/2}} f(x)\|_{L^\infty_{t,|x|}
L^{2+\eta}_\omega}
  \les \|
  e^{i t P^{1/2}} f(x)\|_{L^\infty_t\dot H^{\al}_x} \les \|f\|_{\dot H^{\al}_x}\eneq for $\al\in (1/2, 1]$ and some $\eta>0$.

If we interpolate between \eqref{70-est-WSE-L2} and \eqref{70-est-WSE-Infty} we conclude that,
\beeq\label{70-est-WSE2}\||x|^{\frac{d}{2}-\frac{d+1}{q}-s-\ep} e^{i t P^{1/2}} f(x)
\|_{L^q_{t, |x|\ge R } L^{2+\eta}_\omega }+ \|\<x\>^{-\frac{1}{2}-s-\ep} e^{i t P^{1/2}} f(x)
\|_{L^2_{t, x}}\les
\|f\|_{\dot{H}^{s}_x}\eneq for any $\ep>0$, $2< q\le \infty$ and $s\in (1/2-1/q,1]$ with $\eta >0$ small enough. This implies Theorem \ref{70-thm-WSE} by Duhamel's formula and observation \eqref{70-est-EquiEqn}.

Using the weighted Strichartz estimates, we shall prove the
Strauss conjecture in our setting, by adapting the arguments in
Hidano-Metcalfe-Smith-Sogge-Zhou \cite{HMSSZ}.

To this end, we define $X=X_{s,\ep,q}(\R^d)$ to be the space with norm
defined by
\begin{equation}\label{3.11}
\|h\|_{X_{s,\ep,q}}=\|h\|_{L^{q_s}(|x|\le R)}\, + \, \bigl\| \,
\||x|^{\frac{d}{2}-\frac{d+1}{q}-s-\ep} h\|_{L^q_{|x|}
L^{2+\eta}_\omega(|x|\ge R)} ,\end{equation}if  $
d\bigl(\tfrac12-\tfrac1{q_s})=s$. Note that we have the embedding
$\dot H^{s}\subset X_{s,0,\infty}$ for $s\in(1/2,1]$ by Sobolev
embedding and \eqref{70-est-WSE-Infty}. By duality, we have
\beeq\label{70-est-Embedding}X_{1-s,0,\infty}'\subset \dot H^{s-1}
\textrm{ for }s\in[0,1/2).\eneq We also denote the space
$Y_{s,\ep}(\R^d)$ with norm
$$\|h\|_{Y_{s,\ep}}=\|\<x\>^{-\frac{1}{2}-s-\ep} h
\|_{L^2_{x}}.$$

Note that by Remark \ref{70-rem-Morawetz}, duality and interpolation
and the homogeneous estimate \eqref{70-est-WSE-L2}, we have
\beeq\label{70-est-Morawetz-LowOrder2}\|u\|_{L^2_t
Y_{s,\ep}}+\|u\|_{L^\infty_t \dot{H}^s}+\|\pt u\|_{L^\infty_t
\dot{H}^{s-1}}\les \Vert u_0\Vert_{\dot{H}^{s}}+\Vert u_1\Vert_{\dot
H^{s-1}_x}+ \|F\|_{L^2_t Y_{1-s,\ep}'}\eneq for the solutions to the
linear wave equation \eqref{70-eqn-LW-Stri} and $s\in [0,1]$.

Then by \eqref{70-est-LocEnerDecay} and energy estimate, if
$(\pt^2+P)u=0$, we have
$$\|\phi u\|_{L^p_t \dot{H}^s} \les \|\phi u\|_{L^\infty_t \dot{H}^s}+\|\phi u\|_{L^2_t \dot{H}^s}\les \Vert
u_0\Vert_{\dot{H}^{s}}+\Vert u_1\Vert_{\dot H^{s-1}_x},$$ for any
$\phi\in C_c^\infty$, $s\in[0,1]$ and $p\ge 2$. Thus by
\eqref{70-est-Morawetz-LowOrder2} and Christ-Kiselev lemma, we have
\beeq\label{70-est-LpHs}\|\phi u\|_{L^p_t \dot{H}^s} \les \Vert
u_0\Vert_{\dot{H}^{s}}+\Vert u_1\Vert_{\dot H^{s-1}_x}+\|F\|_{L^2_t
Y_{1-s,\ep}'}\eneq for the solutions of the linear wave equation
\eqref{70-eqn-LW-Stri}, $s\in[0,1]$ and $p> 2$.

In conclusion, by \eqref{70-est-WSE-Infty},
\eqref{70-est-Morawetz-LowOrder2} and \eqref{70-est-LpHs}, we have
\beeq\label{70-est-Strauss-Key}\|u\|_{L^\infty_t \dot H^s\cap L^p_t
X_{s,\ep,p}\cap L^2_t  Y_{s,\ep}}+\|\pt u\|_{L^\infty_t \dot
H^{s-1}}\les \Vert u_0\Vert_{\dot{H}^{s}}+\Vert u_1\Vert_{\dot
H^{s-1}_x}+ \|F\|_{L^2_t Y_{1-s,\ep}'}\eneq for the solutions to the
linear wave equation \eqref{70-eqn-LW-Stri}, if $p>2$, $s\in
(1/2-1/p,1]$. Here, it will be useful to note that if $d=3,4$,
$p>p_c$, i.e. $\frac{d}{2}-\frac{2}{p-1}>\frac{1}{2}-\frac{1}{p}$ and $p>1$, we can choose
$\ep>0$, so that $s=\frac{d}{2}-\frac{2}{p-1}-\frac{p\ep}{p-1}\in (1/2-1/p,1/2)$ and
\beeq\label{70-est-relation}p(\frac{d}{2}-\frac{d+1}{p}-s-\ep)=-(\frac{d}{2}-(1-s))=-s-\frac{d-2}{2}.\eneq

Now we want to prove the higher order estimates for \eqref{70-est-Strauss-Key}. We claim that if
$g^{ij}(x)=h(|x|)\delta^{ij}$,
\begin{eqnarray}
\label{70-est-Strauss-KeyHigh}\sum_{|\al|+|\be|\le 1}\left(
\|\pa_x^\al \Omega^\beta u \|_{L^\infty_t \dot H^s\cap L^p_t
X_{s,\ep,p}\cap L^2_t  Y_{s,\ep}}+\|\pt \pa_x^\al \Omega^\beta
u\|_{L^\infty_t \dot
H^{s-1}}\right) &\les& \\
 \sum_{|\al|+|\be|\le 1} \left(\Vert \pa_x^\al \Omega^\beta u_0\Vert_{\dot{H}^{s}}+
 \Vert \pa_x^\al \Omega^\beta u_1\Vert_{\dot H^{s-1}_x}\right)&&\nonumber\end{eqnarray} for the solutions to the linear wave
equation \eqref{70-eqn-LW-Stri} with $F=0$, if $\rho>0$, $p>2$, $s\in (1/2-1/p,1]$.
As in the proof of \eqref{70-est-Strauss-Key}, we need only to prove
the higher order version of \eqref{70-est-WSE-Infty},
\eqref{70-est-Morawetz-LowOrder2} and \eqref{70-est-LpHs}.

We first estimate $L^2_t  Y_{s,\ep}$ part. By
\eqref{70-est-Morawetz-LowOrder2},  \eqref{a} \eqref{b} in Lemma
\ref{70-lem-equi-PLapBH} and  \eqref{d}  \eqref{e} in Lemma
\ref{70-lem-equi-PLap}, for any $s\in (1/2-1/p,1]$,
\begin{eqnarray*}
\sum_{|\al|\le 1} \|\pa_x^{\al} u\|_{L^2_t  Y_{s,\ep}}&\les& \sum_{|\al|\le 1}
\|\tilde{\pa}_x^{\al} u\|_{L^2_t  Y_{s,\ep}}\\
&\les& \sum_{j\le 1}\|P^{j/2} u\|_{L^2_t  Y_{s,\ep/2}}\\
&\les& \sum_{j\le 1}\left(\Vert P^{j/2} u_0\Vert_{\dot{H}^{s}}+\Vert P^{j/2} u_1\Vert_{\dot H^{s-1}_x}\right)\\
&\les& \sum_{j\le 1}\left(\Vert P^{(j+s)/2} u_0\Vert_{L_x^2}+\Vert P^{(j+s-1)/2} u_1\Vert_{L^2_x}\right)\\
&\les& \sum_{|\al|\le 1}\left(\Vert \pa_x^\al
u_0\Vert_{\dot{H}^{s}}+\Vert \pa_x^\al u_1\Vert_{\dot H^{s-1}_x}\right).
\end{eqnarray*}
Note that since $g^{ij}(x)=h(|x|)\delta^{ij}$,
$[P,\Omega^{ij}]=0$. Thus $$ \sum_{|\al|= 1}\|\Omega^{\al}
u\|_{L^2_t  Y_{s,\ep}}\les  \sum_{|\al| = 1}\left(\Vert \Omega^{\al}
u_0\Vert_{\dot{H}^{s}}+\Vert \Omega^{\al} u_1\Vert_{\dot H^{s-1}_x}\right)$$ In conclusion, we
have for $\rho>0$,
\beeq\label{70-est-Morawetz-LowOrder2Hi}\sum_{|\al|+|\be|\le 1}
\|\pa_x^{\al}\Omega^\be u\|_{L^2_t  Y_{s,\ep}}\les
\sum_{|\al|+|\be|\le 1 }\left(\Vert \pa_x^\al\Omega^\be
u_0\Vert_{\dot{H}^{s}}+\Vert \pa_x^\al\Omega^\be u_1\Vert_{\dot
H^{s-1}_x}\right).\eneq

We turn to the proof of the higher order estimate for
\eqref{70-est-WSE-Infty}. Note that we assume
$g^{ij}(x)=h(|x|)\delta^{ij}$, then $P=-h \Delta+r_1 \pa+r_2$, and
hence by Hardy inequality and \eqref{d} in Lemma
\ref{70-lem-equi-PLap},
$$\|\pa u\|_{\dot{H}^1}\les \|\Delta u\|_{L^2} \les \|P u\|_{L^2} + \|\pa u\|_{L^2} + \|r_2 u\|_{L^2}\les \|P^{1/2} u\|_{H^1}. $$
If we interpolate with \eqref{f}, we get that for $s\in [0,1]$,
\beeq\label{h}\|\pa u\|_{\dot{H}^s}\les \|P^{1/2} u\|_{H^s}.\eneq

 By \eqref{70-est-WSE-Infty}, \eqref{h}, \eqref{d}  and \eqref{e} in Lemma
\ref{70-lem-equi-PLap}, if $(\pt^2+P)u=0$ and $s\in (1/2,1]$, we
have the energy estimate,
\begin{eqnarray*}
\sum_{|\al|\le 1}\|{\pa}^\al_x u\|_{L_t^\infty \dot H^s} &\les&
\sum_{j\le 1}\|P^{j/2} u\|_{L_t^\infty \dot H^s}\\
&\les& \sum_{j\le 1}\left(\|P^{j/2} u(0)\|_{\dot H^s}+\|P^{j/2} \pt
u(0)\|_{\dot H^{s-1}}\right)\\
&\les &\sum_{|\al|\le 1}\left(\| \tilde{\pa}_x^\al u(0)\|_{\dot
H^{s}_x}+\|\tilde {\pa}_x^\al \pt u(0)\|_{\dot H^{s-1}_x}\right)\\
&\les &\sum_{|\al|\le 1}\left(\| {\pa}_x^\al u(0)\|_{\dot
H^{s}_x}+\| {\pa}_x^\al \pt u(0)\|_{\dot H^{s-1}_x}\right)
\end{eqnarray*}
and
\begin{eqnarray*}
\sum_{|\al|\le 1}\||x|^{\frac{d}{2}-s} \tilde{\pa}^\al_x
u\|_{L^\infty_{[0,T],|x|} L^{2+\eta}_\omega}
  &\les& \sum_{|\al|\le 1}\left(\| \tilde{\pa}_x^\al u(0)\|_{\dot
H^{s}_x}+\| \tilde{\pa}_x^\al \pt u(0)\|_{\dot H^{s-1}_x}\right)
+\|[P, \tilde{\pa}_x]u\|_{L^1_T \dot
  H^{s-1}_x}\\
&\les & \sum_{|\al|\le 1}\left(\| {\pa}_x^\al u(0)\|_{\dot
H^{s}_x}+\| {\pa}_x^\al \pt u(0)\|_{\dot
H^{s-1}_x}\right)+\sum_{1\le
|\be|\le 2}\| \tilde{\pa}_x^\be u\|_{L^1_T \dot   H^{s-1}_x}\\
&\les &\sum_{|\al|\le 1}\left(\| {\pa}_x^\al u(0)\|_{\dot
H^{s}_x}+\| {\pa}_x^\al \pt u(0)\|_{\dot H^{s-1}_x}\right)+\sum_{
|\be|\le 1}\| \tilde{\pa}_x^\be u\|_{L^1_T \dot   H^{s}_x}\\
&\les & (1+T)\sum_{|\al|\le 1}\left(\| {\pa}_x^\al u(0)\|_{\dot
H^{s}_x}+\| {\pa}_x^\al \pt u(0)\|_{\dot H^{s-1}_x}\right).
  \end{eqnarray*}
Thus for any $k\in \mathbb{Z}$,
\begin{eqnarray*}\sum_{|\al|\le 1}\||x|^{\frac{d}{2}-s} {\pa}_x^\al
u(t) \|_{L^\infty_{t\in [k, k+1], |x|} L^{2+\eta}_\omega} &\les &
\sum_{|\al|\le 1}\||x|^{\frac{d}{2}-s} \tilde{\pa}^\al_x
u\|_{L^\infty_{t\in[k,k+1],|x|} L^{2+\eta}_\omega}\\
&\les &\sum_{|\al|\le 1}\left(\| {\pa}_x^\al u(k)\|_{\dot
H^{s}_x}+\| {\pa}_x^\al \pt u(k)\|_{\dot H^{s-1}_x}\right)\\
&\les &\sum_{|\al|\le 1}\left(\| {\pa}_x^\al u(0)\|_{\dot
H^{s}_x}+\| {\pa}_x^\al \pt u(0)\|_{\dot H^{s-1}_x}\right).
  \end{eqnarray*}
Hence we get the higher order version of \eqref{70-est-WSE-Infty} if
we combine the commutativity of $P$ and $\Omega$,
\beeq\label{70-est-WSE-InftyHi} \sum_{|\al|+|\be|\le
1}\||x|^{\frac{d}{2}-s} {\pa}_x^\al \Omega^\be u(t) \|_{L^\infty_{t,
|x|} L^{2+\eta}_\omega} \les \sum_{|\al|+|\be|\le 1}\left(\|
{\pa}_x^\al \Omega^\be  u_0\|_{\dot H^{s}_x}+\| {\pa}_x^\al
\Omega^\be u_1\|_{\dot H^{s-1}_x}\right)\eneq for the solution to
the linear wave equation \eqref{70-eqn-LW-Stri} with $F=0$, and $s\in
(1/2,1]$.

We need only to prove the higher order version of \eqref{70-est-LpHs}
now.

By \eqref{70-est-KSS-HighOrder-x}, we know that if $(\pt^2+P)u=0$ and $\rho>0$,
$$\|\<x\>^{-1/2-\ep}\pa_x u\|_{L^2_t \dot H^1}\les \sum_{|\al|= 2}\|\<x\>^{-1/2-\ep}\pa_x^\al u\|_{L^2_{t,x}}+
\|\<x\>^{-3/2-\ep}\pa_x u\|_{L^2_{t,x}}\les \|\pa_x u(0)\|_{H^1_x}+
\| \pt u(0)\|_{H^1_x},$$
$$\|\<x\>^{-1/2-\ep} \pa_x u\|_{L^2_{t,x}}\les \|\pa_x u(0)\|_{L^2_x}+
\|\pt u(0)\|_{L^2_x},$$ Thus for $s\in [0,1]$,
$$\|\phi \pa_x u\|_{L^2_{t}\dot H^s}\les \|\pa
u(0)\|_{H^s_x}+ \|\pt u(0)\|_{H^s_x}\les \sum_{|\al|\le 1} \left(
\|\pa^\al_x u(0)\|_{\dot H^s_x}+ \|\pa_x^\al \pt u(0)\|_{\dot
H^{s-1}_x}\right).$$ Combining the energy estimates, we
know that \beeq\label{70-est-LpHsHi}\sum_{|\al|+|\be|\le 1}\|\phi
\pa^\al_x\Omega^\be u\|_{L^p_t \dot{H}^s} \les \sum_{|\al|+|\be|\le
1}\left(\Vert \pa^\al_x\Omega^\be u_0\Vert_{\dot{H}^{s}}+\Vert
\pa^\al_x\Omega^\be u_1\Vert_{\dot H^{s-1}_x}\right)\eneq for the solutions to the linear
wave equation \eqref{70-eqn-LW-Stri} with $F=0$, $s\in[0,1]$ and $p> 2$. This
completes the proof of \eqref{70-est-Strauss-KeyHigh}.

To conclude this section, let us point out that a similar estimate of \eqref{70-est-Strauss-KeyHigh} holds for the solution to  \eqref{70-eqn-LW}, as in the end of Section \ref{70-sec-Kss-high}.
Precisely, if $g^{ij}(x)=h(|x|)\delta^{ij}$, we have
\begin{eqnarray}
\label{70-est-Strauss-KeyHigh-LB}\sum_{|\al|+|\be|\le 1}\left(
\|\pa_x^\al \Omega^\beta u \|_{L^\infty_t \dot H^s\cap L^p_t
X_{s,\ep,p}\cap L^2_t  Y_{s,\ep}}+\|\pt \pa_x^\al \Omega^\beta
u\|_{L^\infty_t \dot
H^{s-1}}\right) &\les& \\
 \sum_{|\al|+|\be|\le 1} \left(\Vert \pa_x^\al \Omega^\beta u_0\Vert_{\dot{H}^{s}}+
 \Vert \pa_x^\al \Omega^\beta u_1\Vert_{\dot H^{s-1}_x}\right)&&\nonumber\end{eqnarray} for the solutions to the linear wave
equation \eqref{70-eqn-LW} with $F=0$, if $\rho>0$, $p>2$, $s\in (1/2-1/p,1]$.

\section{Strauss Conjecture}
In this section, we prove the Strauss conjecture in the setting where $g^{ij}=h(|x|)\delta^{ij}$, $p>p_c=1+\sqrt{2}$ and $d=3$, i.e., Theorem \ref{70-thm-Strauss}.

By \eqref{70-est-Strauss-KeyHigh-LB}, \eqref{70-est-Embedding} and Duhamel's formula, we have
\begin{eqnarray}
\label{70-est-Strauss-KeyHigh2}&\sum_{|\al|+|\be|\le 1}\left(
\|\pa_x^\al \Omega^\beta u \|_{L^\infty_t \dot H^s\cap L^p_t
X_{s,\ep,p}}+\|\pt \pa_x^\al \Omega^\beta
u\|_{L^\infty_t \dot
H^{s-1}}\right) & \\
\les& \sum_{|\al|+|\be|\le 1} \left(\Vert \pa_x^\al \Omega^\beta u_0\Vert_{\dot{H}^{s}}+
 \Vert \pa_x^\al \Omega^\beta u_1\Vert_{\dot H^{s-1}_x}+
 \|\pa_x^\al \Omega^\beta F\|_{L^1_t \dot H^{s-1}}\right)&\nonumber\\
\les& \sum_{|\al|+|\be|\le 1} \left(\Vert \pa_x^\al \Omega^\beta u_0\Vert_{\dot{H}^{s}}+
 \Vert \pa_x^\al \Omega^\beta u_1\Vert_{\dot H^{s-1}_x}+
 \|\pa_x^\al \Omega^\beta F\|_{L^1_t X_{1-s, 0, \infty}'}\right)&\nonumber
 \end{eqnarray} for the solutions to the linear wave
equation \eqref{70-eqn-LW}, if $\rho>0$, $p>2$, $s\in (1/2-1/p,1/2)$.

Let us now see how we can use these estimates to prove Theorem \ref{70-thm-Strauss}.
We assume  Cauchy data $(u_0,u_1)$ satsifying the smallness condition \eqref{70-eqn-SLW-data},
and let $u^{(0)}$ solve the Cauchy problem \eqref{70-eqn-SLW} with $F=0$.
We iteratively define $u^{(k)}$, for $k\ge 1$, by solving
$$
\begin{cases}
(\partial^2_t-\Delta_{\mathfrak{g}})u^{(k)}(t,x)=F_p(u^{(k-1)}(t,x))\,,
\quad (t,x)\in \R_+\times \Omega
\\
u(0,\cdot)=u_0, \quad
\partial_t u(0,\cdot)=u_1.
\end{cases}
$$
Our aim is to show that if the constant $\de>0$ in \eqref{70-eqn-SLW-data} is
small enough, then so is
$$
M_k = \sum_{|\alpha|\le 1}\left( \|Y^\alpha u^{(k)} \|_{L^\infty_t\dot H^s \cap L^p_t
X_{s,\ep,p}}+ \|\partial_t Y^\alpha u^{(k)}\|_{L^\infty_t \dot H^{s-1}}
\right)
$$for every $k=0,1,2,\dots$

For $k=0$, it follows by \eqref{70-est-Strauss-KeyHigh2} that
$M_0\le C_0\de$, with $C_0$ a fixed constant.
More generally, for some $\eta_1\in (0,1)$, \eqref{70-est-Strauss-KeyHigh2} implies that
\begin{align}\label{3.18}
M_k\le C_0\de +C_0\sum_{|\alpha|\le 1} \,
\Bigl( \,& \bigl\| \, |x|^{-\frac{d}2+1-s}
Y^\alpha F_p(u^{(k-1)})
\bigr\|_{L^1_t L^1_{|x|} L^{2-\eta_1}_\omega(\R_+\times \{|x|\ge R\})}
\\
&+\|Y^\alpha F_p(u^{(k-1)})
\|_{L^1_t L^{q_{1-s}'}_x(\R_+\times \{x\in \Omega: \, |x|\le R\})}\Bigr)\,.
\notag
\end{align}
Note that our assumption \eqref{70-eqn-SLW-Fp} on the nonlinear term $F_p$ implies that
for small $v$
$$
\sum_{|\alpha|\le 1}|Y^\alpha F_p(v)|\lesssim |v|^{p-1}
\sum_{|\alpha|\le 1}|Y^\alpha v|\,.
$$

Since the collection $Y$ contains vectors spanning the tangent space
to $S^{2}$, by Sobolev embedding we have for any $q\le \infty$ and $\eta_2\in(0,1)$
$$
\|v(r\cdot)\|_{L^q_\omega}
\lesssim \sum_{|\alpha|\le 1}
\|Y^\alpha v(r\cdot)\|_{L^{2+\eta_2}_\omega}\,.
$$
Consequently, for fixed $t, r>0$
$$
\sum_{|\alpha|\le 1}\|Y^\alpha F_p(u^{(k-1)}(t,r\cdot) )
\|_{L^{2-\eta_1}_\omega}\lesssim \sum_{|\alpha|\le 1} \|Y^\alpha u^{(k-1)}(t,r\cdot)
\|^p_{L^{2+\eta_2}_\omega}\,.
$$
By \eqref{70-est-relation}, the first
summand in the right side of \eqref{3.18} is dominated by
$C_1 M_{k-1}^p\,.$

We next observe that,
for each fixed $t$, we have
$$
\sum_{|\alpha|\le 1}\|Y^\alpha F_p(u^{(k-1)}(t,\cdot))
\|_{L^{q'_{1-s}}(x:|x|\le R)}
\lesssim \sum_{|\alpha|\le 1} \|u\|^{p-1}_{L^q(x:|x|\le R)} \|Y^\alpha u^{(k-1)}(t,\cdot)
\|_{L^{q_s}(x:|x|\le R)}\,,
$$ where $$\frac{1}{q}=\frac{1}{p-1}(\frac{1}{q_{1-s}'}-\frac{1}{q_s})=\frac{1}{3(p-1)}.$$
It follows by Sobolev embedding on $\{x: |x|\le R\}$ that
$$
\|v\|_{L^q(x :|x|\le R)}
\lesssim \sum_{|\alpha|\le 1}
\|Y^\alpha v\|_{L^{q_s}(x: |x|\le R)}\,,
$$ since $s\in (\frac{1}{2}-\frac{1}{p},\frac{1}{2})\subset [\frac{1}{2}-\frac{1}{p-1}, \frac{3}{2}-\frac{1}{p-1}]$, i.e. $\frac{3}{q}=\frac{1}{p-1}\in [-1+\frac{3}{q_s}, \frac{3}{q_s}]$.

The second summand in the right side of \eqref{3.18} is thus also dominated by
$C_1 M_{k-1}^p$, and we conclude that $M_k\le C_0\de+2 C_0\,C_1 M_{k-1}^p$.
For $\de$ sufficiently small, then
\begin{equation}\label{3.19} M_k\le 2\,C_0\de, \quad k=1,2,3,\dots
\end{equation}

To finish the proof of Theorem \ref{70-thm-Strauss} we need to show that
$u^{(k)}$ converges to a solution of the equation
\eqref{70-eqn-SLW}.  For this it suffices to show that
$$A_k=
\|u^{(k)}-u^{(k-1)}\|_{L^p_t X_{s,\ep,p}}
$$
tends geometrically to zero as $k\to \infty$.  Since
$|F_p(v)-F_p(w)|\lesssim |v-w|(\, |v|^{p-1}+|w|^{p-1}\, )$, the proof of \eqref{3.19} can be adapted to show that,
for small $\de>0$, there is a uniform constant $C$ so that
$$A_k\le CA_{k-1}(M_{k-1}+M_{k-2})^{p-1},$$
which, by \eqref{3.19}, implies that $A_k\le \tfrac12A_{k-1}$ for small
$\de$.  Since $A_1$ is finite, the claim follows, which finishes
the proof of Theorem \ref{70-thm-Strauss}.


\begin{thebibliography}{7}


\bibitem{BaCh99-1}
H.~Bahouri, J.~Y.~Chemin, \textit{Equations d'ondes quasilineaires et estimations de Strichartz (Quasilinear wave
equations and Strichartz estimates)}, Amer. J. Math. \textbf{121} (1999), 1337-–1377.

\bibitem{BaCh99-2}H.~Bahouri, J.~Y.~Chemin, \textit{Equations d'ondes quasilineaires et effet dispersif (Quasilinear
wave equations and dispersive effect)}, Internat. Math. Res. Notices
1999, No. 21, 1141–-1178.

\bibitem{BoHa} J.-F.~Bony, D.~H\"{a}fner,
\textit{The semilinear wave equation on asymptotically Euclidean manifolds}, 
Comm. in PDE \textbf{35}(2010), no. 1, 23-–67.

\bibitem{Burq} N. Burq, \textit{Global Strichartz estimates for nontrapping geometries: about an article by H. F. Smith and C. D. Sogge: ``Global Strichartz estimates for nontrapping perturbations of the Laplacian''}, Comm. Partial Differential Equations, \textbf{28} (2003), 1675--1683.

\bibitem{ChKi01} M.~Christ, A.~Kiselev,
\textit{Maximal Functions Associated to Filtrations}, Jour. Func.
Anal. \textbf{179} (2001), 409--425.



\bibitem{FaWa} D.~Fang, C.~Wang,
\textit{Weighted Strichartz Estimates with Angular Regularity and their
Applications},
Forum Math., to appear.


\bibitem{GLS97} V.~Georgiev, H.~Lindblad, C.~ D.~
Sogge, \textit{Weighted Strichartz estimates and global existence
for semilinear wave equations}, American Journal of Mathematics,
\textbf{119} (1997), 1291--1319.

\bibitem{Glassey} R. T. Glassey,
{\em Existence in the large for $\square u = F(u)$ in two dimensions},
Math. Z. {\bf 178} (1981), 233--261.



\bibitem{HMSSZ} K.~Hidano, J.~Metcalfe, H.~F.~Smith, C.~D.~ Sogge,
Y.~Zhou, \textit{On Abstract Strichartz Estimates and the Strauss
Conjecture for Nontrapping Obstacles},
Trans. Amer. Math. Soc., \textbf{362} (2010), 2789--2809.

\bibitem{Jo81_01}
F.~John, \textit{Blow-up for quasilinear wave equations in three space
  dimensions}, Comm. Pure Appl. Math. \textbf{34} (1981), no.~1, 29--51.

\bibitem{JoKl84_01}
F.~John, S.~Klainerman, \textit{Almost global existence to nonlinear wave
  equations in three space dimensions}, Comm. Pure Appl. Math. \textbf{37}
  (1984), no.~4, 443--455.

\bibitem{JSS}
J.-L. Journ\'e, A. Soffer, C. D. Sogge \textit{Decay estimates for Schr\"odinger
equations} Comm. Pure Appl. Math. \textbf{44} (1991),  573--604.

\bibitem{Ka91}
L.~V.~Kapitanski, \textit{Norm estimates in Besov and Lizorkin-Triebel spaces for the solutions
of second-order linear hyperbolic equations}, J. Soviet Math. \textbf{56} (1991), 2348-–2389.

\bibitem{KeSmSo02} M.~Keel, H.~Smith, C.~D.~Sogge, \textit{Almost global existence for
some semilinear wave equations}, Dedicated to the memory of Thomas
H. Wolff. J. Anal. Math. \textbf{87} (2002), 265--279.

\bibitem{KeSmSo04_01}
M.~Keel, H.~F.~Smith, C.~D.~Sogge, \textit{Almost global existence for quasilinear wave equations in three
  space dimensions}, J. Amer. Math. Soc. \textbf{17} (2004), no.~1, 109--153
  (electronic).


\bibitem{KeTa98} M.~Keel,  T.~Tao,  \textit{Endpoint Strichartz Estimates},
Amer. J. Math. \textbf{120} (1998),  955--980.

\bibitem{KPV} C. E. Kenig, G. Ponce, and L. Vega, {\em On the Zakharov
  and Zakharov-Schulman systems}.  J. Funct. Anal. {\bf 127} (1995), 204--234.

\bibitem{Kl85_01}
S.~Klainerman, \textit{Uniform decay estimates and the {L}orentz invariance of
  the classical wave equation}, Comm. Pure Appl. Math. \textbf{38} (1985),
  no.~3, 321--332.

\bibitem{KlPo83_01}
S.~Klainerman, G.~Ponce, \textit{Global, small amplitude solutions to
  nonlinear evolution equations}, Comm. Pure Appl. Math. \textbf{36} (1983),
  no.~1, 133--141.


\bibitem{LdSo95} H.~Lindblad,  C.~D.~Sogge,  \textit{On Existence and Scattering with
Minimal Regularity for Semilinear Wave Equations},  J. Func. Anal.
\textbf{130} (1995),  357--426.

\bibitem{LS2}{H. Lindblad and C. D. Sogge},
{\em Long-time existence for small amplitude semilinear wave equations},
Amer. J. Math. {\bf 118} (1996), 1047--1135.

\bibitem{Metcalfe} J. Metcalfe, \textit{Global Strichartz estimates for solutions to the wave equation exterior to a convex obstacle},  Trans. Amer. Math. Soc.  \textbf{356}  (2004), 4839--4855.


\bibitem{MeSo06_01}
J.~Metcalfe, C.~Sogge, \textit{Long-time existence of quasilinear wave equations exterior to
  star-shaped obstacles via energy methods}, SIAM J. Math. Anal. \textbf{38}
  (2006), no.~1, 188--209 (electronic).

\bibitem{MeTa07p} J.~Metcalfe, D.~Tataru,
\textit{Global parametrices and dispersive estimates for variable coefficient wave equations}, arXiv:0707.1191.


\bibitem{MoSeSo93}
G.~Mockenhaupt, A.~Seeger, C.~D.~Sogge,
\textit{Local smoothing of {F}ourier integral operators and
  {C}arleson-{S}j\"olin estimates},
J. Amer. Math. Soc., \textbf{6}(1993), no.1, 65--130.


\bibitem{Mor68} C.~S.~Morawetz, \textit{Time decay for the nonlinear Klein-Gordon
equation}, Proc. Roy. Soc. \textbf{A306} (1968), 291--296.

\bibitem{Si83_01}
T.~Sideris, \textit{Global behavior of solutions to nonlinear wave equations in
  three dimensions}, Comm. Partial Differential Equations \textbf{8} (1983),
  no.~12, 1291--1323.

\bibitem{Sideris} T. C. Sideris,
{\em Nonexistence of global solutions to semilinear wave equations
in high dimensions},
J. Differential Equations {\bf 52} (1984), 378--406.

\bibitem{Sm98}
H.~F.~Smith,
\textit{ A parametrix construction for wave equations with {$C\sp {1,1}$}
  coefficients},
 Ann. Inst. Fourier (Grenoble), \textbf{48}(1998), no.3, 797--835.


\bibitem{SmSo00}
H.~F.~Smith, C.~D.~Sogge,  \textit{Global Strichartz estimates for
nontrapping perturbations of the Laplacian}, Comm. Partial
Differential Equations \textbf{25} (2000), no. 11-12, 2171--2183.

\bibitem{Sog95}
C.~D.~Sogge,  \textit{Lectures on nonlinear wave equations. Second edition}.
 International Press,  Boston,  MA,
2008.

\bibitem{StTa} G. Staffilani, D. Tataru, \textit{Strichartz estimates for a Schr\"odinger
equation with nonsmooth coefficients} Comm. Partial Differential Equations \textbf{27} (2002),
1337--1372.

\bibitem{Strauss} W. Strauss, {\em Dispersal of waves vanishing on the
  boundary of an exterior domain}.  Comm. Pure Appl. Math. {\bf 28}
  (1975), 265--278.


\bibitem{Ta00}
D.~Tataru,
\textit{ Strichartz estimates for operators with nonsmooth coefficients and
  the nonlinear wave equation},
Amer. J. Math., \textbf{122}(2000), no.2, 349--376.


\bibitem{Ta01}
D.~Tataru,
\textit{Strichartz estimates for second order hyperbolic operators with
  nonsmooth coefficients. II},
Amer. J. Math., \textbf{123}(2001), no.3, 385--423.

\bibitem{Ta01-2} D.~Tataru, \textit{Strichartz estimates in the hyperbolic space and
global existence for the semilinear wave equation}, Trans. Amer.
Math. Soc. \textbf{353} (2001), no. 2, 795--807.


\bibitem{Ta02} D.~Tataru, \textit{Strichartz estimates for second order hyperbolic operators with nonsmooth coefficients. III}, J. Amer. Math. Soc., \textbf{15}(2):419--442 (electronic), 2002.


\bibitem{Zhou} Y. Zhou,
{\em Cauchy problem for semilinear wave equations with small data in
four space dimensions},
J. Partial Differential Equations {\bf 8} (1995), 135--144.







\end{thebibliography}
\end{document}